\documentclass[leqno,12pt]{article}
\usepackage{latexsym,amsmath,amsfonts,amssymb,mathabx}
\usepackage{enumitem}
\usepackage[mathscr]{eucal}
\usepackage{eqlist}

\usepackage{hyperref}
\hypersetup{
    colorlinks,
    linkcolor={blue},
    citecolor={red},
    urlcolor={blue}
}

\usepackage[usenames,dvipsnames]{color}

\scrollmode

\newtheorem{thm}{Theorem}[section]
\newtheorem{lem}[thm]{Lemma}
\newtheorem{cor}[thm]{Corollary}
\newtheorem{rem}[thm]{Remark}
\newtheorem{rems}[thm]{Remarks}
\newtheorem{prop}[thm]{Proposition}
\newtheorem{df}[thm]{Definition}
\newtheorem{dfs}[thm]{Definitions}
\newtheorem{ex}[thm]{Example}
\newtheorem{exs}[thm]{Examples}
\newtheorem{que}[thm]{Question}
\newtheorem{as}[thm]{Assumption}

\newcommand{\wh}{\widehat}
\newcommand{\wc}{\widecheck}
\newcommand{\wb}{\widebar}

\newcommand{\wt}{\widetilde}
\newcommand{\vS}{\varSigma}
\newcommand{\vO}{\varOmega}
\newcommand{\vT}{\varTheta}
\newcommand{\ov}{\overline}
\newcommand{\mf}{\mathfrak}
\newcommand{\vY}{\varUpsilon}

\newcommand{\E}{\mathbb{E}}
\newcommand{\B}{\mathfrak{B}}

\def\N{{\mathbb N}}
\def\R{{\mathbb R}}
\def\Q{{\mathbb Q}}

\def\F{{\mathcal F}}

\newcommand{\cnp}{\{N_t\}_{t\in\mathbb{R}_{+}}}

\newcommand{\clap}{\{T_n\}_{n\in\N_0}}

\newcommand{\cip}{\{W_n\}_{n\in\N}}

\begin{document}
\title{Some characterizations for Markov processes as mixed renewal processes}
\author{N.D. Macheras and S.M. Tzaninis}

\maketitle

\begin{abstract}
 
In this paper the class of mixed renewal processes (MRPs for short) with mixing parameter  a random vector from \cite{lm6z3} (enlarging Huang's \cite{hu} original class) is replaced by the strictly more comprising class of all extended MRPs by adding a second mixing parameter. We prove under a mild assumption, that within this larger class the basic problem, whether every Markov process is a mixed Poisson process with a random variable as mixing parameter has a solution to the positive. This implies the equivalence of Markov processes, mixed Poisson processes, and processes with the multinomial property within this class. In concrete examples we demonstrate how to establish the Markov property by our results. Another consequence is the invariance of the Markov property under certain changes of measures.\medskip

 \smallskip
\smallskip 

\smallskip

\end{abstract}

\section*{Introduction}\label{intro}

For a given  probability space $(\vO,\vS,P)$ according to Huang \cite{hu}, Definition 3, a mixed renewal process  associated with $\{P_{\wt{y}}\}_{\wt{y}\in\wt\vY}$ and $\nu$ (written $P$-MRP$(\{P_{\wt{y}}\}_{\wt{y}\in\wt\vY},\nu)$ for short), where $\{P_{\wt{y}}\}_{\wt{y}\in\wt\vY}$ is a family of probability measures on $\vS$ and $\nu$ is a probability measure on $\sigma(\{P_{\bullet}(E):E\in\vS\})$, is a counting process $N:=\cnp$ satisfying 
$$ 
P\Bigl(\bigcap_{k=1}^{r}\{W_k\leq w_k\}\Bigr) =\int\prod_{k=1}^{r}P_{\wt{y}}(\{W_k\leq w_k\})\,\nu(d\wt{y}),
$$
if $\{W_n\}_{n\in\N}$ is the interarrival process induced by $N$. In case $\left(P_{\wt{y}}\right)_{W_{n}}={\bf{Exp}}\left(\alpha(\wt{y})\right)$ for some  positive measurable function $\alpha$ on $\R$ a $P$-MRP$(\{P_{\wt{y}}\}_{\wt{y}\in\wt\vY},\nu)$ will be called a mixed Poisson process associated with $\{P_{\wt{y}}\}_{\wt{y}\in\wt\vY}$ and $\nu$ (written $P$-MPP$(\{P_{\wt{y}}\}_{\wt{y}\in\wt\vY},\nu)$ for short).

Under the assumption 
\begin{description}
\item[$(\ast)$]  For $\nu$-almost all $\wt{y}\in\wt\vY$ the function $F_{\wt{y}}:\R_+\longmapsto [0,1]$ defined by means of $F_{\wt{y}}(t):=P_{\wt{y}}(\{W_n\leq t\})$ for all $n\in\N$ is continuously differentiable on $(0,\infty)$ with $0<F'_{\wt{y}}(t)<C$ for each $t>0$, where $C$ is a positive constant, and the function $\alpha:\wt\vY\longmapsto (0,\infty)$ defined by means of $\alpha(\wt y):=\lim_{t\rightarrow 0} F'_{\wt y}(t)$ is measurable,
\end{description}

 Huang (\cite{hu}, Theorem 3) found as his basic result about a $P$-MRP$(\{P_{\wt{y}}\}_{\wt{y}\in\wt\vY},\nu)$ $N$, that $N$ has the Markov property if and only if is a $P$-MPP$(\{P_{\wt{y}}\}_{\wt{y}\in\wt\vY},\nu)$.

An alternative way to model MRPs on  $(\vO,\vS,P)$ within the class of counting processes is to assume the existence of a random vector on the same probability space such that   conditioning on this random vector the counting process behaves like an ordinary renewal process (see \cite{lm6z3}, Definition 3.2 or Definition \ref{rd} (b) of this paper). 

A counting process $N$ being a $P$-MRP$(\{P_{\wt{y}}\}_{\wt{y}\in\wt\vY},\nu)$   is always a MRP according to Definition \ref{rd} (b), while the inverse implication holds true only under additional assumptions (see \cite{lm6z3}, Theorem 4.9).

A special case of a MRP with mixing parameter a random vector is a mixed Poisson process (MPP for short) with mixing parameter a real-valued random variable (written $P$-MPP$(\vT)$ for short) (cf. e.g. \cite{Sc}, page 87 for the definition). It seems that in general there is no relation between a $P$-MPP$(\{P_{\wt{y}}\}_{\wt{y}\in\wt\vY},\nu)$ and a $P$-MPP$(\vT)$.
 
 But under the mild assumption of the existence of a proper disintegration it can be proven that each $P$-MPP$(\vT)$ is a $P$-MPP$(\{P_{\wt{y}}\}_{\wt{y}\in\wt\vY},\nu)$ (see Proposition \ref{prop22}). The inverse implication does not seem to be true without additional assumptions, as it is in general not possible to find for  a given $P$-MPP$(\{P_{\wt{y}}\}_{\wt{y}\in\wt\vY},\nu)$  a real-valued random variable $\vT$ with $P_\vT=\nu$. On the other hand, assuming  there exists such a $\vT$ it is in general not possible to construct conditional probabilities $Q_{\wt y}:=Q(\bullet\mid\vT=\wt{y})$ on $\vS$ such that $Q_{\wt y}=P_{\wt y}$ for $\nu$-a.a. $\wt{y}\in\wt\vY$ with $\wt\vY=R_\vT$, where $R_\vT$ stands for the  range of $\vT$.\smallskip

The above consideration raises the question whether Huang's result can be carried over to MRPs and MPPs with mixing parameter a random vector and a real-valued random variable, respectively.\smallskip

To this purpose, we prove in Section 2  that under a mild assumption a MRP with mixing parameter a random vector is a Markov process if and only if it is a $P$-MPP$(\vT)$ if and only if it has the multinomial property (see Proposition \ref{prop1}).
\smallskip

In Theorem \ref{thm1}, our main result, the above Proposition is generalized for the wider class of extended MRPs (see Definition \ref{rd}, (a)), being strictly more general than the class of Definition \ref{rd}, (b). Proposition \ref{prop1} and Theorem \ref{thm1} are proven under the Assumption \ref{as}, which is essential for the validity of both results (see Remark \ref{r2} (a)).\smallskip

The proofs of Theorem \ref{thm1} and Proposition \ref{prop1} rely on two earlier results of Lyberopoulos and Macheras, where it is proven that under the existence of an appropriate disintegration of $P$ a MRP or a MPP with mixing parameter a random vector or a random variable, respectively, can be reduced to an ordinary renewal  or Poisson process under the disintegrating measures, respectively (see \cite{lm6z3}, Proposition 3.8 and \cite{lm1v}, Proposition 4.4, respectively). Note that the existence of such a disintegration is guaranteed for a wide class of probability spaces used in applied Probability Theory (see the remark following Definition \ref{rcp}). \smallskip

As a consequence of Theorem \ref{thm1}, the invariance of the Markov property, as well as that of the multinomial property, under the change of the measure $P$ into 
disintegrating measures is obtained, see Corollary \ref{cor1}.\smallskip

In Section 3, a method for the construction of non-trivial probability spaces admitting extended MRPs is given, providing concrete examples of probability spaces and extended MRPs satisfying the assumptions of the main result and allowing us to check whether a extended MRP has the Markov property or not.\smallskip

Further applications of our results, concerning the equivalence of the existing definitions of MPPs, are given in the forthcoming paper \cite{lmt2}.

\section{Preliminaries}\label{prel}

By $\N$ is denoted the set of all natural numbers and $\N_0:=\N\cup\{0\}$.  The symbol $\R$ stands for the set of all real numbers, while $\overline{\R}:=\R\cup\{-\infty,+\infty\}$ and $\R^d$ denotes the Euclidean space of dimension $d\in\N$. Given a subset $A$ of a set $\vO$  we denote by $A^c$ the complement $\vO\setminus A$ of $A$ and by $\chi_A$ the indicator function of $A$. For a map $f:D\longmapsto E$ we denote by $R_f$ or by $f(D)$  the set $\{f(x):x\;\in D\}$, and for a set $A\subseteq D$ we denote by $f\upharpoonright A$ the restriction of $f$ to $A$, and by $f(A)$ the set $\{f(x):x\;\in A\}$.

Given a probability space $(\vO,\vS,P)$ a set $N\in\vS$ with $P(N)=0$ is called a $P$-{\bf null set}. For any two sets $A,B\in\vS$ we write $A=_P B$ if $P(A\triangle B)=0$. Given a measurable space $(\vY,H)$, for any two $\vS$-$H$-measurable maps $X,Y:\vO\longmapsto\vY$ we write $X=Y$ $P$-a.s. if $\{X\neq Y\}$ is a $P$-null set.

Given a topology  $\mf{T}$ on $\vO$ write ${\mf B}(\vO)$ for its {\bf Borel $\sigma$-algebra} on $\vO$, i.e. the $\sigma$-algebra generated by $\mf{T}$ and ${\mf B}:={\mf B}(\R)$, ${\overline{\mf B}}:={\mf B}(\overline{\R})$, ${\mf B}_d:={\mf B}(\R^d)$  and $\mf{B}_{\N}:={\mf B}(\R^{\N})$ for the Borel $\sigma$-algebra of subsets of $\R$, $\overline{\R}$, $\R^d$ and $\R^{\N}$ under the corresponding Euclidean topologies, respectively, while ${\mathcal{L}}^{1}(P)$ stands for the family of all real-valued $P$-integrable functions on $\vO$. Functions that are $P$-a.s. equal are not identified.

For the definitions of {\bf real-valued random variables}, {\bf random variables} and {\bf random vectors} we refer to Cohn \cite{Co}, pages 308 and 318.

Given two probability spaces $(\vO,\vS,P)$ and $(\vY,H,Q)$ as well as a $\vS$-$H$-measurable map $X:\vO\longmapsto\vY$ we denote by $\sigma(X):=\{X^{-1}(B): B\in H\}$ the $\sigma$-algebra generated by $X$, while $\sigma(\{X_i\}_{i\in I}):=\sigma\bigl(\bigcup_{i\in{I}}\sigma(X_i)\bigr)$ stands for the $\sigma$-algebra generated by a family $\{X_i\}_{i\in I}$ of $\vS$-$H$-measurable maps from $\vO$ into $\vY$.

For any $d$-dimensional  random vector $X$ on $\vO$ we apply the notation $P_{X}={\bf{K}}(\theta)$ in the meaning that $X$ is distributed according to the law ${\bf{K}}(\theta)$, where $\theta\in\R^d$. In particular, $\mathbf{P}(\theta)$ and $\mathbf{Exp}(\theta)$, where $\theta$ is positive parameter, stand for the law of Poisson and exponential distribution, respectively (cf. e.g. \cite{Sc}).

We write $\E[X|\F]$ for a conditional expectation of $X$ given $\F$ (see \cite{Co}, page 342 for the definition). For $X:=\chi_E\in\mathcal{L}^1(P)$ with $E\in\vS$ we set $P(E\mid\mathcal{F}):=\E_P[\chi_E\mid\mathcal{F}]$.

Given a real-valued random variable $X$ on $\vO$ and a random vector $\vT:\vO\longmapsto\R^d$, a conditional distribution of $X$ over $\vT$ is a map $P_{X|\vT}$ from $\B\times\vO$ into $[0,1]$ such that 
\begin{description}
\item[(cd1)] for each $\omega\in\vO$ the set-function $P_{X|\vT}(\bullet,\omega)$ is a probability measure on $\B$;
\item[(cd2)] for each $B\in\B$ we have
$$
P_{X|\vT}(B,\bullet)=P(\vT^{-1}(B)\mid\sigma(\vT))\quad P\upharpoonright\sigma(\vT)\text{-a.s.,}
$$
where $P_{X|\vT}(B,\bullet)$ is $\sigma(\vT)$-measurable for any fixed $B\in\B$.
\end{description}
For simplicity we write $k:=P_{X|\vT}$ and define the map $K(\vT)$ from $\B\times\vO$ into $[0,1]$ by means of 
$$
K(\vT)(B,\omega):=(k(B,\bullet)\circ\vT)(\omega)\quad \forall\,B\in\B \,\,\,\,\forall\,\omega\in\vO.
$$
Then for $\theta=\vT(\omega)$ with $\omega\in\vO$ the probability measures $k(\bullet,\theta)$ are distributions on $\B$ and so we may write ${\bf K}(\theta)(\bullet)$ instead of $k(\bullet,\theta)$. Consequently, in this case $K(\vT)$ will be written by ${\bf{K}}(\vT)$. 

For any real-valued random variables $X$, $Y$ on $\vO$ we say that $P_{X|\vT}$ and $P_{Y|\vT}$ are $P\upharpoonright\sigma(\vT)$-equivalent and we write $P_{X|\vT}=P_{Y|\vT}$ $P\upharpoonright\sigma(\vT)$-a.s., if there exists a $P$-null set $N\in\sigma(\vT)$ such that for any $\omega\notin N$ and $B\in \B$ the equality $P_{X|\vT}(B,\omega)=P_{Y|\vT}(B,\omega)$ holds true.\smallskip

{\em{From now on $(\vO,\vS,P)$ is a probability space, while $(\vY,H):=(\R,\B)$, $(\varXi,Z):=(\R^d,\B_d)$. Moreover, unless stated otherwise,  $\vT:\vO\longmapsto\R^d$ is a random vector.}}

\section{Characterizations via mixed Poisson processes \\and the multinomial property}\label{mainr} 

We first recall some additional background material, needed in this section.

A family $N:=\cnp$ of  random variables from $(\vO,\vS)$ into $(\ov{\R},\ov{\B})$ is called a {\bf counting process} if there exists a $P$-null set $\vO_N\in\vS$ such that the process $N$ restricted on $\vO\setminus\vO_N$ takes values in $\N_0\cup\{\infty\}$, has right-continuous paths, presents jumps of size (at most) one, vanishes at $t=0$ and increases to infinity. Denote by $T:=\clap$ and $W:=\cip$  the {\bf{arrival process}} and {\bf{interarrival process}} respectively (cf. e.g. \cite{Sc}, Section 1.1, page 6 for the definition) associated with $N$. \smallskip

Recall that a family $\{X_i\}_{i\in I}$ of real-valued random variables  $X_i$ on  $\vO$ 
\begin{itemize}
\item[$\bullet$] is {\bf $P$-conditionally (stochastically) independent given} $\vT$, if for each $n\in\N$ with $n\geq 2$ we have
$$
P(\bigcap^{n}_{j=1} \{X_{i_j}\leq x_{i_j}\}\mid\sigma(\vT))=\prod^{n}_{k=1} P(\{X_{i_j}\leq x_{i_j}\}\mid\sigma(\vT))\qquad\quad P\upharpoonright\sigma(\vT)-\mbox{a.s.}
$$
whenever $i_1,\ldots,i_n$ are distinct members of $I$ and $(x_{i_1},\ldots,x_{i_n})\in\R^n$;
\item[$\bullet$] is {\bf $P$-conditionally identically distributed given} $\vT$, if 
$$
P\bigl(F\cap X_i^{-1}(B)\bigr)=P\bigl(F\cap X_j^{-1}(B)\bigr)
$$
whenever $i,j\in I$, $F\in\sigma(\vT)$ and $B\in \B$.
\end{itemize}

\begin{df}\label{rcp}
\normalfont
Let $Q$ be a probability measure on $\B_d$. A family $\{P_\theta\}_{\theta\in\R^d}$ of probability measures on $\vS$ 
is called a {\bf disintegration} of $P$ over $Q$ if
\begin{description}
\item[(d1)] 
for each $D\in\vS$ the map $\theta\longmapsto P_\theta(D)$ is $\B_d$-measurable;
\item[(d2)]
$\int P_{\theta}(D)\,Q(d\theta)=P(D)$ for each $D\in\vS$.
\end{description}
If $\vT$ is an inverse-measure-preserving function (i.e. $P_\vT(B)=Q(B)$ for each $B\in \B_d$), a disintegration $\{P_{\theta}\}_{\theta\in\R^d}$ of $P$ over $Q$ is called {\bf consistent} with $\vT$ if, for each $B\in \B_d$, the equality $P_{\theta}(\vT^{-1}(B))=1$ holds for $Q$-almost every $\theta\in B$.
\end{df}

{\em Remark.} 
If $\vS$ is countably generated (cf. e.g. \cite{Co}, Section 3.4, page 102 for the definition) and $P$ is perfect (see \cite{fa}, p. 291 for the definition), then there always exists a disintegration $\{P_{\theta}\}_{\theta\in\R^d}$ of $P$ over $Q$ consistent with any inverse-measure-preserving random vector  $\vT:\vO\longmapsto\R^d$ (see \cite{fa}, Theorems 6 and 3).
So, in most cases appearing in applications (e.g. Polish spaces) disintegrations as above always exist.\smallskip

{\em Throughout what follows, unless stated otherwise, $N:=\cnp$ is a counting process, $T:=\left\{T_{n}\right\}_{n\in\N_0}$ is an  arrival process, $W:=\left\{W_{n}\right\}_{n\in\N}$ is its induced  interarrival process  and without loss of generality we may and do assume that $\vO_N=\emptyset$. Moreover, we simply write ``conditionally'' in the place of ``conditionally given $\vT$'' whenever $\vT$ is clear from the context.}\smallskip

A Poisson process $N$ under $P$ with parameter $\theta>0$ is denoted by $P$-PP$(\theta)$.

The following result shows that under a mild assumption every $P$-MPP($\vT$) is a $P$-MPP$(\{P_{\wt{y}}\}_{\wt{y}\in\wt\vY},\nu)$. But as pointed out in the introduction the inverse implication is in general not possible.

\begin{prop}\label{prop22}
Let $\wb{\vT}$ be a real valued random variable such that $N$ is a $P$-MPP$(\wb\vT)$, and suppose that $\{R_{\wb\theta}\}_{\wb\theta\in\R}$ is a disintegration of $P$ over $P_{\wb\vT}$ consistent with $\wb\vT$. Then $N$ is a $P$-MPP$(\{R_{\wb\theta}\}_{\wb\theta\in\R},P_{\wb\vT})$.
\end{prop}

{\bf Proof.}  For every $r\in\N$ and for all $w_{1},\ldots ,w_{r}\in(0,\infty)$ we get  
\begin{eqnarray*}
P\biggl(\bigcap^{r}_{k=1}\left\{W_{k}\leq w_{k}\right\}\biggr)&=&\int_{\vO} P\left(\bigcap^{r}_{k=1}\left\{W_{k}\leq w_{k}\right\}\mid \wb{\vT}\right)dP\\
&=&\int_{(0,\infty)} R_{\wb{\theta}}\left(\bigcap^{r}_{k=1}\{W_k\leq w_k\}\right)P_{\wb{\vT}}(d\wb{\theta})\\
&=&\int_{(0,\infty)}\prod^{r}_{k=1}\,(1- e^{-\wb{\theta} w_k})\,P_{\wb{\vT}}(d\wb{\theta}),
\end{eqnarray*}
where the second equality follows from \cite{lm1v}, Lemma 3.5, and the third one follows from \cite{lm1v}, Proposition 4.4, together with \cite{Sc}, Theorem 2.3.4. \hfill$\Box$\smallskip

It seems to be natural to generalize  the notion of a MPP$(\vT)$  through the next definition (a), just by adding to the structural parameter $\vT$ another ``mixing"  parameter $h$. Then the resulting {\bf extended MRPs} also comprise the processes studied in \cite{lm1v}.

\begin{dfs}\label{rd}
\normalfont
{\bf{(a)}} A counting process $N$ is called an  {\bf extended MRP  with mixing parameters $\vT$ and $h$, and interarrival time conditional distribution ${\bf K}(h(\vT))$} (written $P$-eMRP$({\bf K}(h(\vT)))$ for short), if $h$ is a $\R^k$-valued $\mathfrak B(D)$-$\mathfrak B_k$-measurable function on $D\in\B_d$ with $R_\vT\subseteq D$ for $k\in\N$, if the induced  interarrival process $W$ is $P$-conditionally independent and 
$$
\forall\;n\in\N\qquad [P_{W_n|\vT}={\bf{K}}\left(h(\vT)\right)\qquad P\upharpoonright\sigma(\vT)-\text{a.s.}].
$$
Without loss of generality we may and do assume that 
\begin{equation}\label{rd1}
\forall\, n\in\N\qquad [P_{W_n|\vT}={\bf{K}}\left(h(\vT)\right)].
\end{equation}

{\bf{(b)}} An extended MRP  with mixing parameters $\vT$ and $\mbox{id}_D$ is called a {\bf MRP with mixing parameter $\vT$} (written $P$-MRP$({\bf{K}}(\vT))$ for short).

In particular, if there exists a $\theta_0\in\R^d$ with $P(\{\vT=\theta_0\})=1$, then $N$ is  a {\em renewal process} with  interarrival time distribution $\mathbf{K}(\theta_0)$ (written $P$-RP$(\mathbf{K}(\theta_0))$ for short).
\end{dfs}

Huang's Definition 3 from \cite{hu}  at that time did not involve a structural parameter $\vT$ and as a result it is strictly less general than that of a MRP$({\bf{K}}(\vT))$ as witnessed by Theorem 4.9 from \cite{lm6z3}.\smallskip

\begin{it}
If no confusion arises, we may write eMRP$({\bf K}(h(\vT)))$, MRP$({\bf{K}}(\vT))$,  RP$({\bf{K}}(\theta_0))$, \\
MRP$(\{P_{\wt{y}}\},\nu)$, MPP$(\vT)$, PP$(\theta)$ and MPP$(\{P_{\wt{y}}\},\nu)$ in the place of $P$-eMRP$({\bf K}(h(\vT)))$, $P$-MRP$({\bf{K}}(\vT))$, $P$-RP$({\bf{K}}(\theta_0))$, $P$-MRP$(\{P_{\wt{y}}\}_{\wt{y}\in\wt\vY},\nu)$, $P$-MPP$(\vT)$, $P$-PP$(\theta)$ and \\
$P$-MPP$(\{P_{\wt{y}}\}_{\wt{y}\in\wt\vY},\nu)$ respectively.\smallskip
\end{it}

{\em From now on, unless stated otherwise, $\{P_\theta\}_{\theta\in D}$ is a disintegration of $P$ over $P_\vT$ consistent with $\vT$, where $D\in\B_d$.}\smallskip

Before we  formulate the basic result of this section we need the next auxiliary lemma.

\begin{lem}\label{lem1}
Let $D\in\B_d$ with $R_\vT\subseteq D$ and let $h: D\longmapsto \R^k\;(k\in\N)$ be a $\mathfrak B(D)$-$\mathfrak B_k$-measurable function. Assume that there exists a $P_{\vT}$-null set $L_0\in \mathfrak B(D)$ such that the restriction $h\upharpoonright D\setminus L_0$ is injective. Put $\wt{\vT}:=h\circ \vT$, $g:=\left(h\upharpoonright D\setminus L_0\right)^{-1}: h\left(D\setminus L_0\right)\longmapsto D\setminus L_0$, and $M:=h(D\setminus L_0)$. For any $\wt{\theta}\in\R^k$ and $A\in\vS$ define
$$
Q_{\wt{\theta}}\left(A\right):=\begin{cases}(P_{\bullet}(A)\circ g)(\wt{\theta}) &\text{if}\;\; \wt{\theta}\in M; \\ P\left(A\right)& \text{if}\;\; \wt{\theta}\in \R^k\setminus M.\quad\quad\end{cases}
$$
Then the following holds true:
\begin{enumerate}
\item 
the family $\left\{Q_{\wt{\theta}}\right\}_{\wt{\theta}\in\R^k}$ is a disintegration of $P$ over $P_{\wt{\vT}}$ consistent with $\wt{\vT}$;
\item
for any $n\in\N$ the equivalence
$$
\forall\;\; \theta\in D\setminus L_0 \quad [\left(P_{\theta}\right)_{W_{n}}={\bf{K}}\left(h(\theta)\right)]\Longleftrightarrow \forall\;\; \wt{\theta}\in M\quad [\left(Q_{\wt{\theta}}\right)_{W_{n}}={\bf{K}}(\wt{\theta})]
$$
is fulfilled;
\item
for every $\theta\in D\setminus L_0$ the process $W$ is $P_{\theta}$-independent if and only if for every $\wt{\theta}\in M$ it is $Q_{\wt{\theta}}$-independent.
\end{enumerate}
\end{lem}
{\bf{Proof.}} First note that $M\in\B_k$ (cf. e.g. \cite{Co}  Theorem 8.3.7), the function  $g$ is $\mathfrak B(M)$-$\mathfrak B(D\setminus L_0)$-measurable (cf. e.g. \cite{Co} Proposition 8.3.5), and $P_{\wt{\vT}}(M)=1$.

\textit{Ad (i):} Clearly, $\left\{Q_{\wt{\theta}}\right\}_{\wt{\theta}\in\R^k}$ is a family of probability measure on $\vS$ satisfying condition (d1). Condition (d2) follows by (d2) for $\left\{P_{\theta}\right\}_{\theta\in D}$.  

To show that $\left\{Q_{\wt{\theta}}\right\}_{\wt{\theta}\in\R^k}$ is consistent with $\wt{\vT}$, let $A\in\vS$ and $B\in\mathfrak B_k$ be arbitrary. Putting $E:=h^{-1}\left(B\right)$ we have
\begin{equation}
B\cap M=_{P_{\wt\vT}} g^{-1}(E\cap (D\setminus L_0)).
\label{eq:2}
\end{equation}

Thus,

\begin{eqnarray*}
\int_B Q_{\wt{\theta}}(A)P_{\wt{\vT}}(d\wt{\theta})&=&\int_{B\cap M} Q_{\wt{\theta}}(A)P_{\wt{\vT}}(d\wt{\theta})+\int_{B\cap (\R^k\setminus M)} Q_{\wt{\theta}}(A)P_{\wt{\vT}}(d\wt{\theta})\\
&\stackrel{(\ref{eq:2})}{=}&\int_{g^{-1}(E\cap (D\setminus L_0))} Q_{\wt{\theta}}(A)P_{\wt{\vT}}(d\wt{\theta})\\
&=&\int_{g^{-1}\left(E\cap (D\setminus L_0)\right)} \left(P_{\bullet}(A)\circ g\right)(\wt{\theta})P_{\wt{\vT}}(d\wt{\theta})\\
&=&\int_{E\cap (D\setminus L_0)} P_{\theta}(A)P_{\vT}\left(d\theta\right)=P\left(A\cap \vT^{-1}\left(E\cap (D\setminus L_0)\right)\right)\\
&=&P\left(A\cap (\wt{\vT})^{-1}\left(B\right)\right),
\end{eqnarray*}

where the sixth equality follows by the consistency of $\left\{P_{\theta}\right\}_{\theta\in D}$ with $\vT$.
This completes the proof of \textit{(i)}.\smallskip

\textit{Ad (ii):} Let us fix on arbitrary $n\in\N$ and $A\in\vS$ such that $A:=W_n^{-1}(B)$ for  $B\in\mathfrak B\left((0,\infty)\right)$. Assume that for all $\theta\in  D\setminus L_0$ we have $\left(P_{\theta}\right)_{W_{n}}={\bf{K}}\left(h(\theta)\right)$. Then for any $\wt{\theta}\in M$ we get $\theta:=g(\wt{\theta})\in  D\setminus L_0$, implying that 
\begin{eqnarray*}
\left(Q_{\wt{\theta}}\right)_{W_n}(B)&=&Q_{\wt{\theta}}(A)=\left(P_{\bullet}(A)\circ g\right)(\wt{\theta})= P_{g(\wt{\theta})}(A)=P_{\theta}(A)\\
&=& \left(P_{\theta}\right)_{W_n}(B)= {\bf{K}}\left(h(\theta)\right)(B)={\bf{K}}(\wt{\theta})(B).
\end{eqnarray*}

For the inverse implication, assume that for any $\wt{\theta}\in M$ we have $Q_{\wt{\theta}}(A)={\bf{K}}(\wt{\theta})$. Then for any $\theta\in D\setminus L_0$ we get $\wt{\theta}:=g^{-1}(\theta)\in M$; hence 
\begin{eqnarray*}
\left(P_{\theta}\right)_{W_n}(B)&=& 
P_{\theta}(A)= P_{g(\wt{\theta})}(A)=\left(P_{\bullet}(A)\circ g\right)(\wt{\theta})= Q_{\wt{\theta}}(A)\\
&=&\left(Q_{\wt{\theta}}\right)_{W_n}(B)={\bf{K}}(\wt{\theta})(B)={\bf{K}}\left(h(\theta)\right)(B).
\end{eqnarray*}
Assertion \textit{(iii)} follows in a similar way. 
\hfill$\Box$

\begin{rem}
\label{rem1}
\normalfont
The following result is well known (cf. e.g. \cite{hu}, Theorem 2) but we write it exactly in the form, that we need.\smallskip

\begin{it}
 Let $\theta\in\R^d$ be fixed and let $N$ be a RP$(\mathbf{K}(\theta))$. For any $t\in\R_+$ put $F_{\theta}(t):= P(\{W_{n}\leq t\})$ for all $n\in\mathbb N$. Assume that the function $F_{\theta}$ is continuously differentiable on $(0,\infty)$, $0<F'_{\theta}(t)<C$ for each $t>0$, where $C$ is a positive constant, which may depend of $\theta$, and that $p_d\left(\theta\right):=\lim_{t\rightarrow 0}F'_{\theta}(t)$ is positive. 
Then the following assertions are equivalent:
\begin{enumerate}
\item
$N$ has the Markov property;
\item
$N$ is a PP$(p_d(\theta))$. 
\end{enumerate}
\end{it}
\end{rem}

It is well known that if $N$ is a MPP$(\vT)$ then it satisfies the Markov property (cf. e.g. \cite{Sc}, Theorem 4.2.3 and page 44 for the definition of the Markov property). However, the trivial counting process $N$ defined by means of $N_{t}:= \left[t\right]$ for every $t\in\R_+$, where by $\left[t\right]$ is denoted the integer part of $t$,  is a Markov RP$(\mathbf{K}(\theta_0))$ but not any Poisson process. 
 This raises the question, under which conditions  a Markov MRP$({\bf{K}}(\vT))$ is a MPP$(\vT)$ ? \smallskip
 
 Under the following mild assumption this question is answered to the positive in Proposition \ref{prop1}.
 
\begin{as}\label{as}
\normalfont

Let $D\in\B$ with $R_\vT\subseteq D$, $h:D\longmapsto\R$ be a $\B(D)$-measurable function, let $N$ be a $P$-eMRP$({\bf K}(h(\vT)))$ and let $\{P_\theta\}_{\theta\in D}$  be  a disintegration of $P$ over $P_{\vT}$ consistent with $\vT$. It follows by \cite{lm6z3}, Lemma 3.5 together with condition (\ref{rd1}) that 
\marginpar{rd2}
\begin{equation}\label{rd2}
\forall\, n\in\N\quad\forall\,\theta\in{D}\qquad [(P_{\theta})_{W_n}={\bf{K}}\left(h(\theta)\right)].
\end{equation}
For any $\theta\in D$ and $t\in\R_+$ put
$$
F_{h(\theta)}(t):=P_\theta(\{W_n\leq t\})\quad\text{for all}\quad n\in\N.
$$
Clearly the function $F_{h(\theta)}$ depends on the distribution of $W_n$ and, because of condition \eqref{rd2}, on $h$. We say that $N$, $h$ and $\{P_\theta\}_{\theta\in D}$  satisfy Assumption \ref{as}, if there exists a $P_\vT$-null set $L_h:=L_{h,N,\{P_\theta\}_{\theta\in D}}$ in $\mathfrak B(D)$ such that for any $\theta\notin L_{h}$ the function $F_{h(\theta)}$ is continuously differentiable on $(0,\infty)$, there exists a function $C\in\mathcal L^1(P_{h(\vT)})$ with $0<F'_{h(\theta)}(t)<C(h(\theta))$ for each $t>0$, and the function $p_h:D\setminus L_h\longmapsto\R$ defined by means of $p_h(\theta):=p_{h,1}(\theta):=\lim_{t\rightarrow 0} F'_{h(\theta)}(t)$ is positive and injective.

For the special case  $D=\R$ and $h:=id_{\R}$ we write for simplicity $L$, $F_\theta$ and $p_1$ in the place of $L_h$, $F_{h(\theta)}$ and $p_{h}$ respectively, and we say that $N$ and $\{P_\theta\}_{\theta\in\R}$  satisfy Assumption \ref{as}.
\end{as}


\begin{prop}\label{prop1}
Consider the following statements.
\begin{enumerate}
\item
$N$ has the multinomial property;
\item
$N$ has the Markov property;
\item 
$N$ is a MPP$(\widecheck{\vT})$.
\end{enumerate}

Then $(iii)\Rightarrow(i)\Rightarrow(ii)$.
If $N$ is a MRP$({\bf{K}}(\vT))$ and $\{P_{\theta}\}_{\theta\in\R^d}$ is  a disintegration of $P$ over $P_{\Theta}$ consistent with $\vT$ satisfying Assumption \ref{as}, put $\wc\vT(\omega):=(p_d\circ\vT)(\omega)$ if $\omega\in\vT^{-1}(L^c)$, and denote again by $\wc\vT$ any measurable extension of $\wc\vT$ from $\vT^{-1}(L^c)$ to $\vO$, then (i) to (iii) are all equivalent.
\end{prop}

{\bf Proof.} For the definition of the multinomial property cf. e.g. \cite{SZ}, page 2 or \cite{Sc}, Lemma 2.3.1. The implication $(i)\Longrightarrow (ii)$ follows by an easy computation, while the implication $(iii)\Longrightarrow{(i)}$ is well known for any real-valued random variable on $\vO$ in the place of $\widecheck{\vT}$ (cf. e.g. \cite{Sc}, Lemma 4.2.2).\smallskip

Ad $(ii)\Longrightarrow (iii)$: Let $N$ be a MRP$({\bf K}(\vT))$ having the Markov property, such that $N$ and $\{P_\theta\}_{\theta\in\R^d}$ satisfy Assumption \ref{as}, and let $\wc\vT$ be as above.
\smallskip

{\bf{(a)}} There exists a $P_{\vT}$-null set $L_1\in \mathfrak B_d$ such that $N$ is a $P_\theta$-RP$({\bf{K}}(\theta))$ for all $\theta\notin L_1$.\smallskip

In fact, since $N$ is a MRP$({\bf{K}}(\vT))$ we may apply \cite{lm6z3}, Proposition 3.8, to obtain (a).\smallskip

{\bf{(b)}} For every $t>0$ and $n\in\mathbb N_{0}$ condition $P\left(\{N_{t}=n\}\right)>0$ holds true.
\smallskip

In fact, first notice that by applying \cite{lm1v}, Lemma 3.5, for any $n\in\mathbb N_0$ and any $t>0$  we obtain 
$$P(\{N_{t}=n\})=\mathbb E_{P_{\vT}}\left[P_{\bullet}(\{N_{t}=n\})\right].$$
Therefore, it is sufficient to show that  for every $t>0$ and $n\in\mathbb N_{0}$ condition  $P_{\theta}(\{N_{t}=n\})>0$ holds true for $P_{\vT}$-a.a. $\theta\in\R^d$. But this can be easily shown by induction on $n$.
\smallskip

{\bf{({c})}} For every $s,t\in (0,\infty)$ with $s<t$ condition
$$P\left(\{N_s=N_t=1\}\right)=-\mathbb E_{P_\vT}\left[\int_{0}^{s}G_{\bullet}(t-x)G'_{\bullet}(x)dx\right]>0,$$
where $G_{\theta}(t):=1-F_{\theta}(t)$ for any  $\theta\in \R^d$, 
 holds true.\smallskip

In fact, for every $s,t\in (0,\infty)$ with $s<t$  
we obtain
\begin{eqnarray*}
P(\{N_{s}=N_{t}=1\})&=&P(\{T_{1}\leq s<t<T_{2}\})\\
&=&P(\{W_{1}\leq s\})-P(\{W_{1}\leq s,W_{2}\leq t-W_{1}\})\\
&=&\int_{\R^d\cap(L\cup L_1)^{c}}\left[\int^{s}_{0}f_{\theta}(x)dx-\int^{s}_{0}\int^{t-x}_{0}f_{\theta}(y)f_{\theta}(x)dydx\right]P_{\vT}(d\theta)\\
&=&-\mathbb E_{P_{\vT}}\left[\int^{s}_{0}G_{\bullet}(t-x)G'_{\bullet}(x)dx\right],
\end{eqnarray*}
where the third equality follows from 
\cite{lm1v}, Lemma 3.5, and  (a), and $f_\theta:=F'_\theta$. Assume now, if possible, that $P(\{N_{s}=N_{t}=1\})=0$. 
Due to  Assumption \ref{as}, the latter is equivalent to the fact that for all $\theta\notin L\cup L_1$ and $x\in(0,s]$ condition 
$G_{\theta}(t-x)= 0$ holds true; hence  $f_{\theta}(t-x)= 0$ for every $x\in(0,s]$, a contradiction according  to Assumption \ref{as}. \smallskip

Note that, due to  (b) and (c) all conditional probabilities considered in the next three steps are well defined.

The proofs of the following  steps (d), (e) and (g) consist of some modifications of the corresponding arguments of Huang in the proof of Theorem 3 from \cite{hu}. We include the detailed proofs for the sake of completeness.
\smallskip

{\bf{(d)}} For any $0<t,v$ condition 
\begin{equation}
\frac{\mathbb E_{P_{\vT}}\left[-G_{\bullet}(t+v)p_d(\bullet)\right]}{\mathbb E_{P_{\vT}}\left[G_{\bullet}(v)G'_{\bullet}(t)\right]}=\frac{\mathbb E_{P_{\vT}}\left[-G_{\bullet}(t)p_d(\bullet)\right]}{\mathbb E_{P_{\vT}}\left[G'_{\bullet}(t)\right]}=1,
\label{eq:h2}
\end{equation}
holds true.\smallskip

In fact, for any $0<u<t$ and $v>0$ applying the Markov property we have 
\begin{eqnarray*}
\lefteqn{P(\{N_{t-u}=N_{t}=N_{t+v}=1\})=P(\{N_{t+v}=1\}|\{N_{t}=1\})P(\{N_{t-u}=N_{t}=1\})}&&\\
&\Longleftrightarrow &\frac{P(\{N_{t-u}=N_{t}=N_{t+v}=1\})}{P(\{N_{t-u}=N_{t}=1\})}=P(\{N_{t+v}=1\}|\{N_{t}=1\})\\
&\Longleftrightarrow &\frac{P(\{N_{t-u}=N_{t+v}=1\})}{P(\{N_{t-u}=N_{t}=1\})}=\frac{A_{1}(t,v)}{B_{1}(t)}
\end{eqnarray*}
where $A_{1}(t,v):=P(\{N_{t}=N_{t+v}=1\})$ and $B_{1}(t):=P(\{N_{t}=1\})$. Then, by  (c) we obtain that
$$P(\{N_{t-u}=N_{t+v}=1\})=-\mathbb E_{P_{\vT}}\left[\int^{t-u}_{0}G_{\bullet}(t+v-x)G'_{\bullet}(x)dx\right]$$
and 
$$P(\{N_{t-u}=N_{t}=1\})=-\mathbb E_{P_{\vT}}\left[\int^{t-u}_{0}G_{\bullet}(t-x)G'_{\bullet}(x)dx\right].$$

Thus,
$$
\frac{\mathbb E_{P_{\vT}}\left[\int^{t-u}_{0}G_{\bullet}(t+v-x)G'_{\bullet}(x)dx\right]}{\mathbb E_{P_{\vT}}\left[\int^{t-u}_{0}G_{\bullet}(t-x)G'_{\bullet}(x)dx\right]}=\frac{A_{1}(t,v)}{B_{1}(t)}.
$$

It is obvious that the right side of the last equation is independent of $u$. Therefore its derivative with respect to  $u$ must be equal to zero; hence, equating the derivative of the left side with zero and  applying the Dominated Convergence Theorem we have
\begin{equation*}
\frac{\mathbb E_{P_{\vT}}\Big[\int^{t-u}_{0}G_{\bullet}(t+v-x)G'_{\bullet}(x)dx\Big]}{\mathbb E_{P_{\vT}}\Big[\int^{t-u}_{0}G_{\bullet}(t-x)G'_{\bullet}(x)dx\Big]}=\frac{\mathbb E_{P_{\vT}}\Big[G_{\bullet}(u+v)G'_{\bullet}(t-u)\Big]}{\mathbb E_{P_{\vT}}\Big[G_{\bullet}(u)G'_{\bullet}(t-u)\Big]}.
\end{equation*}
Since the left side of the last equality is independent of $u$ the same must hold for the right one. Thus, letting on the right side $u\rightarrow 0$ and $u\rightarrow t$ by the Dominated Convergence Theorem  we obtain
\begin{equation}
\frac{\mathbb E_{P_{\vT}}\Big[-G_{\bullet}(t+v)p_d(\bullet)\Big]}{\mathbb E_{P_{\vT}}\Big[G_{\bullet}(v)G'_{\bullet}(t)\Big]}=\frac{\mathbb E_{P_{\vT}}\Big[-G_{\bullet}(t)p_d(\bullet)\Big]}{\mathbb E_{P_{\vT}}\Big[G'_{\bullet}(t)\Big]}.
\label{eq:h3}
\end{equation}
Moreover, since the right side of the (\ref{eq:h3}) is independent of $v$ its derivative with respect to $v$ must be equal to zero. Thus, equating the derivative of the left side  with zero and applying the Dominated Convergence Theorem we have
\begin{equation*}
\frac{\mathbb E_{P_{\vT}}\left[-G_{\bullet}(t+v)p_d(\bullet)\right]}{\mathbb E_{P_{\vT}}\left[G_{\bullet}(v)G'_{\bullet}(t)\right]}=\frac{\mathbb E_{P_{\vT}}\left[-G'_{\bullet}(t+v)p_d(\bullet)\right]}{\mathbb E_{P_{\vT}}\left[G'_{\bullet}(v)G'_{\bullet}(t)\right]}.
\end{equation*}
Thus, for $v\rightarrow 0$ and by the Dominated Convergence Theorem we obtain that
\begin{equation*}
\frac{\mathbb E_{P_{\vT}}\left[-G_{\bullet}(t)p_d(\bullet)\right]}{\mathbb E_{P_{\vT}}\left[G'_{\bullet}(t)\right]}=\frac{\mathbb E_{P_{\vT}}\left[-G'_{\bullet}(t)p_d(\bullet)\right]}{\mathbb E_{P_{\vT}}\left[-p_d(\bullet)G'_{\bullet}(t)\right]}=1;
\end{equation*}
implying together with condition (\ref{eq:h3}) that
\begin{equation*}
\frac{\mathbb E_{P_{\vT}}\Big[-G_{\bullet}(t+v)p_d(\bullet)\Big]}{\mathbb E_{P_{\vT}}\Big[G_{\bullet}(v)G'_{\bullet}(t)\Big]}=\frac{\mathbb E_{P_{\vT}}\Big[-G_{\bullet}(t)p_d(\bullet)\Big]}{\mathbb E_{P_{\vT}}\Big[G'_{\bullet}(t)\Big]}=1.
\end{equation*}

{\bf{(e)}} For any $0<t,v$ condition 
\begin{equation}
\mathbb E_{P_{\vT}}\left[G_{\bullet}(t+v)\left(p_d(\bullet)\right)^{2}\right]=\mathbb E_{P_{\vT}}\left[-G_{\bullet}(v)G'_{\bullet}(t)p_d(\bullet)\right]
\label{eq:h4}
\end{equation}
holds true.\smallskip

In fact, for any $0<u<t$, $v>0$ and $w>0$  applying the Markov property we have 
\begin{eqnarray*}
\lefteqn{P(\{N_{t-u}=N_{t}=1,N_{t+v}=N_{t+v+w}=2\})}\\&&=P(\{N_{t+v+w}=2\}|\{N_{t+v}=2\})P(\{N_{t+v}=2\}|\{N_{t}=1\})P(\{N_{t-u}=N_{t}=1\})
\end{eqnarray*}
or equivalently
\begin{equation*}
\frac{P(\{N_{t-u}=N_{t}=1,N_{t+v}=N_{t+v+w}=2)\}}{P(\{N_{t-u}=N_{t}=1)\}}=\frac{A_{2}(t,v,w)}{B_{2}(t,v)},
\end{equation*}

where $A_{2}(t,v,w)=P(\{N_{t+v+w}=2,N_{t+v}=2\})\cdot P(\{N_{t+v}=2,N_{t}=1\})$ and $B_{2}(t,v)=P(\{N_{t+v}=2\})\cdot P(\{N_{t}=1\})$.

Moreover working as in the proof of (c) we get 
\begin{eqnarray*}
\lefteqn{P(\{N_{t-u}=N_{t}=1,N_{t+v}=N_{t+v+w}=2\})}\\
&&=\mathbb E_{P_{\vT}}\left[\int^{t-u}_{0}\int^{t+v-x}_{t-x}G_{\bullet}(t+v+w-x-y)G'_{\bullet}(y)G'_{\bullet}(x)dydx\right],
\end{eqnarray*}
implying together with (c) that 
\begin{equation*}
\frac{\mathbb E_{P_{\vT}}\left[\int^{t-u}_{0}\int^{t+v-x}_{t-x}G_{\bullet}(t+v+w-x-y)G'_{\bullet}(y)G'_{\bullet}(x)dydx\right]}{\mathbb E_{P_{\vT}}\left[\int^{t-u}_{0}G_{\bullet}(t-x)G'_{\bullet}(x)dx\right]}=-\frac{A_{2}(t,v,w)}{B_{2}(t,v)}.
\end{equation*}
It is obvious that the right side of the last equation is independent of $u$. Therefore its derivative with respect to  $u$ must be equal to zero. Furthermore, equating the left's side derivative with zero and  applying the Dominated Convergence Theorem we have
\begin{eqnarray*}
\lefteqn{\frac{\mathbb E_{P_{\vT}}\left[\int^{t-u}_{0}\int^{t+v-x}_{t-x}G_{\bullet}(t+v+w-x-y)G'_{\bullet}(y)G'_{\bullet}(x)dydx\right]}{\mathbb E_{P_{\vT}}\left[\int^{t-u}_{0}G_{\bullet}(t-x)G'_{\bullet}(x)dx\right]}}\\
&&\quad\quad\quad=\frac{\mathbb E_{P_{\vT}}\left[\int^{u+v}_{u}G_{\bullet}(u+v+w-y)G'_{\bullet}(y)G'_{\bullet}(t-u)dy\right]}{\mathbb E_{P_{\vT}}\left[G_{\bullet}(u)G'_{\bullet}(t-u)\right]}.
\end{eqnarray*}
Since the left side of the above equality is independent of $u$ the same must hold for the right one; hence, taking on the right side $u\rightarrow 0$ and $u\rightarrow t$ by the Dominated Convergence Theorem we obtain
\begin{eqnarray*}
\frac{\mathbb E_{P_{\vT}}\left[-\int^{t+v}_{t}G_{\bullet}(t+v+w-y)G'_{\bullet}(y)p_d(\bullet)dy\right]}{\mathbb E_{P_{\vT}}\left[\int^{v}_{0}G_{\bullet}(v+w-y)G'_{\bullet}(y)G'_{\bullet}(t)dy\right]}=\frac{\mathbb E_{P_{\vT}}\left[-G_{\bullet}(t)p_d(\bullet)\right]}{\mathbb E_{P_{\vT}}\left[G'_{\bullet}(t)\right]}=1,
\end{eqnarray*}
where the last equality can be rewritten in the following form
\begin{equation*}
\mathbb E_{P_{\vT}}\left[-\int^{t+v}_{t}G_{\bullet}(t+v+w-y)G'_{\bullet}(y)p_d(\bullet)dy\right]=\mathbb E_{P_{\vT}}\left[\int^{v}_{0}G_{\bullet}(v+w-y)G'_{\bullet}(y)G'_{\bullet}(t)dy\right]
\end{equation*}
due  to  (d). 
If we take the derivative with respect to $v$ and we apply the Dominated Convergence Theorem in the above equality we obtain that 
\begin{equation*}
\mathbb E_{P_{\vT}}\left[-G_{\bullet}(w)G'_{\bullet}(t+v)p_d(\bullet)\right]=\mathbb E_{P_{\vT}}\left[G_{\bullet}(w)G'_{\bullet}(v)G'_{\bullet}(t)\right].
\end{equation*}
Derivating now with respect to $w$ and applying once again the Dominated Convergence Theorem we get that 
\begin{equation*}
\mathbb E_{P_{\vT}}\left[-G'_{\bullet}(w)G'_{\bullet}(t+v)p_d(\bullet)\right]=\mathbb E_{P_{\vT}}\left[G'_{\bullet}(w)G'_{\bullet}(v)G'_{\bullet}(t)\right].
\end{equation*}
By letting now $w\rightarrow 0$ and by the Dominated Convergence Theorem we obtain that 
\begin{equation*}
\mathbb E_{P_{\vT}}\left[G'_{\bullet}(t+v)\left(p_d(\bullet)\right)^{2}\right]=\mathbb E_{P_{\vT}}\left[-G'_{\bullet}(v)G'_{\bullet}(t)p_d(\bullet)\right].
\end{equation*}
Finally, integration with respect to $v$ yields
\begin{equation*}
\mathbb E_{P_{\vT}}\left[G_{\bullet}(t+v)\left(p_d(\bullet)\right)^{2}\right]=\mathbb E_{P_{\vT}}\left[-G_{\bullet}(v)G'_{\bullet}(t)p_d(\bullet)\right].
\end{equation*}

{\bf{(f)}} For any $0<t,v$ condition 
\begin{equation}
\mathbb E_{P_{\vT}}\left[G_{\bullet}(t)G_{\bullet}(v)\left(p_d(\bullet)\right)^{2}\right]=\mathbb E_{P_{\vT}}\left[-G_{\bullet}(v)G'_{\bullet}(t)p_d(\bullet)\right]
\label{eq:h6}
\end{equation}
holds true.\smallskip

In fact, for any $0<u<t$, $v>0$ and $w>0$ by the Markov property we have 
\begin{eqnarray*}
\lefteqn{P(\{N_{t-u}=N_{t}=1,N_{t+v}=2,N_{t+v+w}=3\})}\\&&=P(\{N_{t+v+w}=3\}|\{N_{t+v}=2\})P(\{N_{t+v}=2\}|\{N_{t}=1\})P(\{N_{t-u}=N_{t}=1\}),
\end{eqnarray*}
or equivalently
\begin{equation*}
\frac{P(\{N_{t-u}=N_{t}=1,N_{t+v}=2,N_{t+v+w}=3\})}{P(\{N_{t-u}=N_{t}=1\})}=\frac{A_{3}(t,v,w)}{B_{3}(t,v)},
\end{equation*}
where $A_{3}(t,v,w)=P(\{N_{t+v+w}=3,N_{t+v}=2\})\cdot P(\{N_{t+v}=2,N_{t}=1\})$ and $B_{3}(t,v)=P(\{N_{t+v}=2\})\cdot P(\{N_{t}=1\})$.

Moreover, applying  Lemma 3.5, of \cite{lm1v} and  (a) after some manipulation as in the proof of (c) we obtain
\begin{eqnarray*}
\lefteqn{P(\{N_{t-u}=N_{t}=1,N_{t+v}=2,N_{t+v+w}=3\})}\\
&=&\mathbb E_{P_{\vT}}\Big[\int^{t-u}_{0}\int^{t+v-x}_{t-x}\int^{t+v+w-x-y}_{t+v-x-y}G_{\bullet}(t+v+w-x-y-z) f_{\bullet}(z)f_{\bullet}(y)f_{\bullet}(x)dzdydx\Big].
\end{eqnarray*}
The latter together with  (c) yields
\begin{eqnarray*}
\lefteqn{\frac{\mathbb E_{P_{\vT}}\left[\int^{t-u}_{0}\int^{t+v-x}_{t-x}\int^{t+v+w-x-y}_{t+v-x-y}G_{\bullet}(t+v+w-x-y-z)G'_{\bullet}(z)G'_{\bullet}(y)G'_{\bullet}(x)dzdydx\right]}{\mathbb E_{P_{\vT}}\left[\int^{t-u}_{0}G_{\bullet}(t-x)G'_{\bullet}(x)dx\right]}}\\
&&\quad\quad\quad\quad\quad\quad\quad\quad\quad\quad\quad\quad\quad\quad\quad\quad\quad\quad\quad\quad\quad\quad\quad\quad\quad\quad\quad\quad=\frac{A_{3}(t,v,w)}{B_{3}(t,v)}.
\end{eqnarray*}
It is obvious that the right side of the last equation is independent of $u$. Therefore its derivative with respect to  $u$ must be equal to zero; hence, equating the derivative of the left side  with zero and  applying the Dominated Convergence Theorem we have
\begin{eqnarray*}
\lefteqn{\frac{\mathbb E_{P_{\vT}}\left[\int^{t-u}_{0}\int^{t+v-x}_{t-x}\int^{t+v+w-x-y}_{t+v-x-y}G_{\bullet}(t+v+w-x-y-z)f_{\bullet}(z)f_{\bullet}(y)f_{\bullet}(x)dzdydx\right]}{\mathbb E_{P_{\vT}}\left[\int^{t-u}_{0}G_{\bullet}(t-x)G'_{\bullet}(x)dx\right]}}\\
&&=\frac{\mathbb E_{P_{\vT}}\left[\int^{u+v}_{u}\int^{u+v+w-y}_{u+v-y}G_{\bullet}(u+v+w-y-z)G'_{\bullet}(z)G'_{\bullet}(y)G'_{\bullet}(t-u)dzdy\right]}{\mathbb E_{P_{\vT}}\left[G_{\bullet}(u)G'_{\bullet}(t-u)\right]}.
\end{eqnarray*}
Since the left side of the last equality is independent of $u$ the same must hold for the right one. Consequently, taking on the right side  $u\rightarrow 0$ and $u\rightarrow t$ by the Dominated Convergence Theorem we obtain
\begin{eqnarray*}
\frac{\mathbb E_{P_{\vT}}\left[-\int^{t+v}_{t}\int^{t+v+w-y}_{t+v-y}G_{\bullet}(t+v+w-y-z)G'_{\bullet}(z)G'_{\bullet}(y)p_d(\bullet)dzdy\right]}{\mathbb E_{P_{\vT}}\left[\int^{v}_{0}\int^{v+w-y}_{v-y}G_{\bullet}(v+w-y-z)G'_{\bullet}(z)G'_{\bullet}(y)G'_{\bullet}(t)dzdy\right]}=\frac{\mathbb E_{P_{\vT}}\left[-G_{\bullet}(t)p_d(\bullet)\right]}{\mathbb E_{P_{\vT}}\left[G'_{\bullet}(t)\right]}
\end{eqnarray*}
Due to condition (\ref{eq:h2}), the last equality can be rewritten in the following form
\begin{eqnarray*}
\lefteqn{\mathbb E_{P_{\vT}}\left[\int^{v}_{0}\int^{v+w-y}_{v-y}G_{\bullet}(v+w-y-z)G'_{\bullet}(z)G'_{\bullet}(y)G'_{\bullet}(t)dzdy\right]}\\
&&=\mathbb E_{P_{\vT}}\left[-\int^{t+v}_{t}\int^{t+v+w-y}_{t+v-y}G_{\bullet}(t+v+w-y-z)G'_{\bullet}(z)G'_{\bullet}(y)p_d(\bullet)dzdy\right].
\end{eqnarray*}
If we take the derivative with respect to $v$ and we apply the Dominated Convergence Theorem in the above equality we obtain that 
\begin{eqnarray*}
\lefteqn{\mathbb E_{P_{\vT}}\left[\int^{w}_{0}G_{\bullet}(w-z)G'_{\bullet}(z)G'_{\bullet}(v)G'_{\bullet}(t)dz\right]}\\
&&=\mathbb E_{P_{\vT}}\left[-\int^{w}_{0}G_{\bullet}(w-z)G'_{\bullet}(z)G'_{\bullet}(t+v)p_d(\bullet)dz\right],
\end{eqnarray*}
or equivalently if we put $s:=w-z$ in the first integral we obtain
\begin{eqnarray*}
\lefteqn{\mathbb E_{P_{\vT}}\left[\int^{w}_{0}G_{\bullet}(s)G'_{\bullet}(w-s)G'_{\bullet}(v)G'_{\bullet}(t)ds\right]}\\
&&=\mathbb E_{P_{\vT}}\left[-\int^{w}_{0}G_{\bullet}(w-z)G'_{\bullet}(z)G'_{\bullet}(t+v)p_d(\bullet)dz\right].
\end{eqnarray*}
Derivating now with respect to $w$ and applying once again the Dominated Convergence Theorem we get that 
\begin{equation*}
\mathbb E_{P_{\vT}}\left[-G_{\bullet}(w)p_d(\bullet)G'_{\bullet}(v)G'_{\bullet}(t)\right]=\mathbb E_{P_{\vT}}\left[-G'_{\bullet}(w)p_d(\bullet)G'_{\bullet}(t+v)\right].
\end{equation*}
Integration with respect to $t$ yields
\begin{equation*}
\mathbb E_{P_{\vT}}\left[-G_{\bullet}(t)G_{\bullet}(w)p_d(\bullet)G'_{\bullet}(v)\right]=\mathbb E_{P_{\vT}}\left[-G_{\bullet}(t+v)G'_{\bullet}(w)p_d(\bullet)\right].
\end{equation*}
For $v\rightarrow 0$ by the Dominated Convergence Theorem we get 
\begin{equation*}
\mathbb E_{P_{\vT}}\left[G_{\bullet}(t)G_{\bullet}(w)\left(p_d(\bullet)\right)^{2}\right]=\mathbb E_{P_{\vT}}\left[-G_{\bullet}(t)G'_{\bullet}(w)p_d(\bullet)\right].
\end{equation*}
Take now $v\rightarrow t$, the by the Dominated Convergence Theorem we get that
\begin{equation*}
\mathbb E_{P_{\vT}}\left[G_{\bullet}(v)G_{\bullet}(w)\left(p_d(\bullet)\right)^{2}\right]=\mathbb E_{P_{\vT}}\left[-G_{\bullet}(v)G'_{\bullet}(w)p_d(\bullet)\right].
\end{equation*}
Finally by letting $w\rightarrow t$ and by the Dominated Convergence Theorem  we obtain
\begin{equation*}
\mathbb E_{P_{\vT}}\left[G_{\bullet}(t)G_{\bullet}(v)\left(p_d(\bullet)\right)^{2}\right]=\mathbb E_{P_{\vT}}\left[-G_{\bullet}(v)G'_{\bullet}(t)p_d(\bullet)\right].
\end{equation*}

{\bf{(g)}} For every $t>0$ condition  
\begin{equation}
\mathbb E_{P_{\vT}}\left[\big(G'_{\bullet}(t)+p_d(\bullet)G_{\bullet}(t)\big)^{2}\right]=0
\label{eq:h8}
\end{equation}
holds true.\smallskip

In fact, let us fix on arbitrary $t>0$. Applying (\ref{eq:h4}) and (\ref{eq:h6}) for $v\rightarrow t$ we get 
\begin{eqnarray*}
\lefteqn{\mathbb E_{P_{\vT}}\left[\big(G'_{\bullet}(t)+  p_d(\bullet)G_{\bullet}(t)\big)^{2}\right]}\\
&&=\mathbb E_{P_{\vT}}\left[G'^{2}_{\bullet}(t)-G^{2}_{\bullet}(t)p_d(\bullet)^2-G_{\bullet}(2t)p_d(\bullet)^2+ p_d(\bullet)^2 G^{2}_{\bullet}(t)\right];
\end{eqnarray*}
hence
\begin{equation}
\mathbb E_{P_{\vT}}\left[\big(G'_{\bullet}(t)+p_d(\bullet)G_{\bullet}(t)\big)^{2}\right]=\mathbb E_{P_{\vT}}\left[G'^{2}_{\bullet}(t)-G_{\bullet}(2t)p_d(\bullet)^2\right].
\label{eq:h9}
\end{equation}
Derivating  (\ref{eq:h2}) with respect to $v$ and applying the Dominated Convergence Theorem we get
\begin{equation*}
\mathbb E_{P_{\vT}}\left[-G'_{\bullet}(t+v)p_d(\bullet)\right]=\mathbb E_{P_{\vT}}\left[G'_{\bullet}(v)G'_{\bullet}(t)\right];
\end{equation*}
hence by letting $v\rightarrow t$ we obtain 
\begin{equation*}
\mathbb E_{P_{\vT}}\left[-G'_{\bullet}(2t)p_d(\bullet)\right]=\mathbb E_{P_{\vT}}\left[G'^{2}_{\bullet}(t)\right],
\end{equation*}
implying that equation (\ref{eq:h9})  can be rewritten as
\begin{equation}
\mathbb E_{P_{\vT}}\left[\big(G'_{\bullet}(t)+p_d(\bullet)G_{\bullet}(t)\big)^{2}\right]=\mathbb E_{P_{\vT}}\left[-G'_{\bullet}(2t)p_d(\bullet)-G_{\bullet}(2t)p_d(\bullet)^2\right].
\label{eq:h10}
\end{equation}
Taking now $v\rightarrow 0$ in equation (\ref{eq:h4}) we obtain that 
\begin{equation*}
\mathbb E_{P_{\vT}}\left[G_{\bullet}(t)p_d(\bullet)^2\right]=\mathbb E_{P_{\vT}}\left[-G'_{\bullet}(t)p_d(\bullet)\right];
\end{equation*}
 hence substituting $t$ by $2t$ we get that 
\begin{equation}
\mathbb E_{P_{\vT}}\left[G_{\bullet}(2t)p_d(\bullet)^2\right]=\mathbb E_{P_{\vT}}\left[-G'_{\bullet}(2t)p_d(\bullet)\right].
\label{eq:h11}
\end{equation}
Thus, by equations (\ref{eq:h10}) and (\ref{eq:h11}) it follows that 
\begin{equation*}
\mathbb E_{P_{\vT}}\left[\big(G'_{\bullet}(t)+p_d(\bullet)G_{\bullet}(t)\big)^{2}\right]=0\quad\text{for any}\;\;t>0.
\end{equation*}

{\bf{(h)}} There exists a $P_{\vT}$-null set $L_2\in \mathfrak B_d$, containing $L$, such that for any $\theta\notin L_2$ the process $W$ is $P_{\theta}$-exponentially distributed with parameter $p_d\left(\theta\right)$.\smallskip

In fact, step (g) yields that for any $s\in\Q_+\setminus\{0\}$ 
there exists a $P_{\vT}$-null set $M_s\in \mathfrak B_d$, containing $L$, such that for any $\theta\notin M_s$ condition 
\begin{equation}
G'_{\theta}(s)+p_d(\theta)G_{\theta}(s)=0
\label{eq:ha}
\end{equation}
holds true. Put $L_2:=\bigcup_{s\in\Q_+\setminus\{0\}} M_s$ and let $t>0$ and $\theta\notin L_2$ be arbitrary. There exists a sequence $\{s_n\}_{n\in\N}$ in $\Q_+\setminus\{0\}$ such that $t=\lim_{n\rightarrow \infty} s_n$. Applying (\ref{eq:ha}) we get $G'_{\theta}(s_n)=-p_d(\theta)G_{\theta}(s_n)$, implying together with Assumption \ref{as} that  $G'_{\theta}(t)=\lim_{n\rightarrow \infty}G'_{\theta}(s_n)=-p_d(\theta)G_{\theta}(t)$; hence  $G_{\theta}(t)=e^{-p_d(\theta)t}$, or $\left(P_{\theta}\right)_{W_n}={\bf{Exp}}\left(p_d(\theta)\right)$ for all $n\in\N$. 
\smallskip

{\bf{(i)}} Put $L_\ast:=L_1\cup L_2$ and denote again by $p_d$ the restriction of $p_d$ to $L_\ast^c$. Define $M_\ast:=p_d(L_\ast)$ and $r:=p_d^{-1}:p_d\left(L_\ast^c\right)\longmapsto L_\ast^c$, as well as the family $\{Q_{\widecheck{\theta}}\}_{\widecheck{\theta}\in\R}$ as in Lemma  \ref{lem1} (for $k=1$, and $p_d$, $r$ and $L_\ast$ in place of $h$, $g$ and $L_0$, respectively). Then the family $\{Q_{\widecheck{\theta}}\}_{\widecheck{\theta}\in\R}$ is a disintegration of $P$ over $P_{\widecheck{\vT}}$ consistent with $\widecheck{\vT}$,  and for any $\widecheck{\theta}\notin M_\ast$  the process $W$ is $Q_{\widecheck{\theta}}$-independent and $(Q_{\widecheck{\theta}})_{W_n}={\bf Exp}(\widecheck{\theta})$ for any $n\in\N$. \smallskip

In fact, by Lemma \ref{lem1} the family $\{Q_{\widecheck{\theta}}\}_{\widecheck{\theta}\in\R}$ is a disintegration of $P$ over $P_{\widecheck{\vT}}$ consistent with $\widecheck{\vT}$ and for any $\widecheck{\theta}\notin M_\ast$ taking into account  (a) and (h)  
we get that $W$ is $Q_{\widecheck{\theta}}$-independent and $(Q_{\widecheck{\theta}})_{W_n}={\bf Exp}(\widecheck{\theta})$ for any $n\in\N$.\smallskip

{\bf{(j)}} $N$ is a MPP$(\widecheck{\vT})$.\smallskip

In fact, by  (i) we obtain that for any $\widecheck{\theta}\notin M_\ast$ the counting process $N$ is $Q_{\widecheck{\theta}}$-PP$(\widecheck{\theta})$  (cf. e.g \cite{Sc}, Theorem 2.3.4). Thus, applying \cite{lm1v}, Proposition 4.4, we deduce that $N$ is a MPP$(\widecheck{\vT})$.
\hfill$\Box$

\begin{rems}\label{r1}
\normalfont
{\bf (a)} Assumption \ref{as} is a modification of Huang's Assumption $(\ast)$ since there it is  assumed the stronger condition $0<F'_\theta(t)<C$ for any $t>0$ and $\theta\notin L$, where $C$ is a positive constant, in the place of  $0<F'_\theta(t)<C(h(\theta))$ for any $t>0$ and $\theta\notin L$ for $C\in\mathcal L^1(P_{h(\vT)})$ in Assumption \ref{as}. In addition in Assumption \ref{as} it is assumed that $p_h$ is measurable and almost everywhere  positive and injective function, while Huang assumes that $p_h$ is a measurable and positive function.

We do not know if our assumption about the injectivity of $p_h$ is essential. We need the injectivity assumption for the proofs of Proposition \ref{prop1} and Theorem \ref{thm1}. These rely on Lemma \ref{lem1} and it would be interesting for us to know whether the injectivity assumption of $p_h$ in Lemma \ref{lem1} is essential.\smallskip

{\bf{(b)}} Our proof of step (f) of Proposition \ref{prop1} differs from that of  Huang's \cite{hu} proof for condition (16) of Theorem 3, since we found it impossible to actually carry through the suggestions given by Huang for this step.

In fact, starting from the probability $P(\{N_{t-u}=N_{t}=1,N_{t+v}=N_{t+v+w}=2\})$, as suggested by Huang, we could prove only condition (15) instead of condition (16) of Huang \cite{hu} as it is shown in step (e) of the proof of Proposition \ref{prop1}.\smallskip

{\bf (c)} The arguments used in the proof of the implication $(ii)\Longrightarrow (iii)$ of Proposition \ref{prop1} are totally different from those used in the proof of \cite{SZ}, Theorem 4.2, since we do not use the Bernstein-Widder Theorem.\smallskip

{\bf (d)} Huang's Theorem 3 in \cite{hu} remains true under the weaker assumption $0<F'_{\wt y}(t)<C(h(\wt y))$ for any $t>0$ and $\nu$-a.a. $\wt y\in \wt\vY$, where $C\in\mathcal L^1(\nu)$, in the place of $0<F'_{\wt\vY}(t)<C$ for any $t>0$ and $\nu$-a.a. $\wt y\in \wt\vY$, where $C$ is a positive constant of Huang's assumption $(\ast)$.
\end{rems}

The following result extends Lemma 3.5 from \cite{lm6z3}.

\begin{lem}\label{lem2}
Let $D\in\B_d$ with  $R_\vT\subseteq D$, let $h: D\longmapsto \R^k$ be a $\B(D)$-$\B_k$-measurable function, and let $\{X_n\}_{n\in\N}$ be a sequence of real-valued random variables. Then the following are equivalent 
\begin{enumerate}
\item
$\exists\;\;L_3\in\sigma\left(\vT\right)_0\;\;\forall\;\;n\in\N\;\;[P_{X_n|\vT}\upharpoonright \B\times L_3^c={\bf{K}}_n\left(h(\vT)\right)\upharpoonright \B\times L_3^c]$;
\item
$\exists\;\;\wt{L}_3\in \mathfrak B(D)_0\;\;\forall\;\;n\in\N\;\;\forall\;\;\theta\in D\setminus\wt{L}_3\;\;[\left(P_\theta\right)_{X_n}={\bf{K}}_n\left(h(\theta)\right)],$
\end{enumerate} 
where $\sigma\left(\vT\right)_0:=\{M\in\sigma\left(\vT\right):P(M)=0\}$ and $\mathfrak B(D)_0:=\{\wt{M}\in \mathfrak B(D): P_{\vT}(\wt{M})=0\}$.
\end{lem}

{\bf{Proof.}} Ad \textit{(i)}$\Longrightarrow$\textit{(ii)}: Assume that there exists a set $L_3\in\sigma\left(\vT\right)_0$ such that for every $n\in\N$ condition 
$$
P_{X_n|\vT}\upharpoonright \B\times L_3^c={\bf{K}}_n\left(h(\vT)\right)\upharpoonright \B\times L_3^c
$$
holds true. Then for any fixed $n\in\N$, $F\in \mathfrak B_k$ and $B\in \B$ we obtain
$$
\int_{\vT^{-1}(F)} P_{X_n|\vT}(B,\bullet)dP=\int_{\vT^{-1}(F)} {\bf{K}}_n\left(h(\vT)\right)(B,\bullet)dP;
$$
hence taking into account \cite{lm1v}, Lemma 3.5, we get 
$$
\int_{F} \left(P_{\theta}\right)_{X_n}(B)P_{\vT}(d\theta)=\int_{F} {\bf{K}}_n(h(\theta))(B)P_{\vT}(d\theta)
$$
Consequently, there exists a set $\wt{L}_{n,B}\in \mathfrak B(D)_0$ such that 
\begin{equation}
\left(P_{\theta}\right)_{X_n}(B)=  {\bf{K}}_n(h(\theta))(B)\qquad\text{for any}\;\;\theta\in D\setminus\wt{L}_{n,B}.
\label{eq:lem21}
\end{equation}
Put $\wt{L}_3:=\bigcup_{n\in\N}\bigcup_{B\in\mathcal G_\B} \wt{L}_{n,B}$, where $\mathcal G_\B$ is a countable generator of $\B$ being closed under finite intersections, and denote by $\mathcal D$ the class of all $B\in \B$ such that condition (\ref{eq:lem21}) is satisfied for every $\theta\in D\setminus\wt{L}_3$ and $n\in\N$. It can be easily seen that $\mathcal G_\B\subseteq \mathcal D$ and that $\mathcal D$ is a Dynkin class, implying that $\mathcal D=\B$. Thus assertion \textit{(ii)} follows.

Applying a similar reasoning we obtain the converse implication.
\hfill$\Box$\smallskip

The following result shows how to reduce a eMRP$({\bf K}(h(\vT)))$ to a MRP$({\bf{K}}(\wt{\vT}))$ under the change of the mixing parameter.

\begin{lem}\label{lem3}
Let $h$  and $\wt{\vT}$ be as in Lemma \ref{lem1}. Suppose that $N$ is a eMRP$({\bf K}(h(\vT)))$ on $(\vO,\vS,P)$. Then $N$ is a MRP$({\bf{K}}(\wt{\vT}))$.
\end{lem}

{\bf{Proof.}} Let $\{P_{\theta}\}_{\theta\in D}$, $g$ and $\{Q_{\wt{\theta}}\}_{\wt{\theta}\in\R^k}$ be as in Lemma \ref{lem1}. According to Lemma \ref{lem2}, there exists a $P_\vT$-null set $\wt{L}_3\in \mathfrak B(D)$  such that $\left(P_{\theta}\right)_{W_n}={\bf{K}}\left(h(\theta)\right)$ for all $\theta\in D\setminus\wt{L}_3$.
We may and do assume that $\wt{L}_3$ contains the $P_\vT$-null set $L_{0}$ of Lemma \ref{lem1}. Applying now Lemma \ref{lem1}, we obtain that $\{Q_{\wt{\theta}}\}_{\wt{\theta}\in\R^k}$ is a disintegration of $P$ over $P_{\wt{\vT}}$ consistent with $\wt{\vT}$, and for all $\wt{\theta}\in h(D\setminus\wt{L}_3)$ the process $W$ is $Q_{\wt{\theta}}$-independent and $\left(Q_{\wt{\theta}}\right)_{W_n}={\bf{K}}(\wt{\theta})$ for every $n\in\N$. But the latter together with \cite{lm6z3} Proposition 3.8 yields the conclusion of the lemma.
\hfill$\Box$\smallskip

The next result extends Proposition \ref{prop1}.

\begin{thm}\label{thm1}
Let $h$, $L_0$ be  as in Lemma \ref{lem1}, let $N$ be an eMRP$({\bf K}(h(\vT)))$ on $(\vO,\vS,P)$ satisfying together with $\{P_{\theta}\}_{\theta\in D}$ Assumption \ref{as}. Put $O_h:=L_0\cup L_h$ and $\wh{\vT}(\omega)=(p_h\circ\vT)(\omega)$ if $\omega\in\vT^{-1}(D\setminus O_h)$, and denote again by $\wh\vT$ any measurable extension of $\wh\vT$ from  $\vT^{-1}(D\setminus O_h)$ to $\vO$. Then the following are all equivalent:
\begin{enumerate}
\item
$N$ has the multinomial property;
\item
$N$ has the Markov property;
\item
$N$ is a MPP$(\wh{\vT})$. 
\end{enumerate}
\end{thm}

{\bf{Proof.}} Let $\wt{\vT}$ and $\{Q_{\wt{\theta}}\}_{\wt{\theta}\in\R^k}$ be as in Lemma \ref{lem1}. It then follows by Lemma \ref{lem3} that the process $N$ is a MRP$({\bf{K}}(\wt{\vT}))$. For any $\wt{\theta}\in\R^k$ and $t\in\R_+$ put $F_{\wt{\theta}}(t):=Q_{\wt{\theta}}(\{W_n\leq t\})$ for all $n\in\N$.\smallskip

{\bf{(a)}} Put $V_h:=h(D\setminus O_h)$. Then $P_{\wt\vT}(V_h)=1$ and for any $\wt{\theta}\in V_h$ the function $F_{\wt{\theta}}$ is continuously differentiable on $(0,\infty)$ and $0<F'_{\wt{\theta}}(t)<C(\wt{\theta})$ for each $t>0$.
\smallskip

In fact, according to Assumption \ref{as} there exists a $P_{\vT}$-null set $L_h\in \mathfrak B(D)$ such that for each $\theta\in D\setminus L_h$ the function $F_{h(\theta)}$ is continuously differentiable on $(0,\infty)$ and $0<F'_{h(\theta)}(t)<C(h(\theta))$, implying that the same holds true for each $\theta\in D\setminus O_h$. By Lemma \ref{lem1} for any $\wt{\theta}\in V_h$ condition $\left(Q_{\wt{\theta}}\right)_{W_n}={\bf{K}}(\wt{\theta})$ holds true for any $n\in\N$. Then for any $\wt{\theta}\in V_h$ there exists exactly one $\theta\in D\setminus O_h$ such that $\wt{\theta}=h(\theta)$ and  
\begin{equation}
F_{\wt{\theta}}(t)={\bf{K}}(\wt{\theta})\left((-\infty,t]\right)={\bf{K}}\left(h(\theta)\right)\left((-\infty,t]\right)=F_{h(\theta)}(t)
\label{eq:thm11}
\end{equation}
for any $t\in\R_+$, implying that  (a) holds true.\smallskip

{\bf{(b)}} The function $p_k:V_h\longmapsto\R$ defined by $p_k(\wt{\theta}):=\lim_{t\rightarrow 0} F'_{\wt{\theta}}(t)$ for any $\wt{\theta}\in V_h$ is positive and injective.\smallskip

In fact, since $N$, $h$ and $\{P_\theta\}_{\theta\in D}$ satisfy Assumption \ref{as}, taking into account condition (\ref{eq:thm11}) we get that $p_k(\wt{\theta})=\lim_{t\rightarrow 0} F'_{h(\theta)}(t)=p_h(\theta)>0.$

Clearly, $\wh{\vT}(\omega)=(p_h\circ\vT)(\omega)=(p_k\circ h\circ \vT)(\omega)=(p_k\circ\wt{\vT})(\omega)$ for any $\omega\in\vT^{-1}(D\setminus O_h)$, and due to (a) and (b) we deduce that $N$ and $\{Q_{\wt{\theta}}\}_{\wt{\theta}\in\R^k}$ satisfy Assumption \ref{as}; hence by Proposition \ref{prop1} for $\wt{\vT}$ and $p_k$ in the place of $\vT$ and $p_d$, respectively, we get the thesis of the theorem. 
\hfill$\Box$\smallskip

The following result may be of independent interest, since it ensures the permanence of the Markov and the multinomial property with respect to $P$ to that with respect to the disintegrating measures $Q_{\wh{\theta}}$.

\begin{cor}\label{cor1}
Let $\{P_{\theta}\}_{\theta\in D}$, $h$, $N$, $F_{h(\theta)}$, $p_h$ and $\wh{\vT}$  be as in Theorem \ref{thm1}. Fix on an arbitrary $A\in\vS$ and for any $\wh{\theta}\in\R$  put
$$
Q_{\wh{\theta}}\left(A\right):=\begin{cases}(P_{\bullet}(A)\circ p^{-1}_h)(\wh{\theta}) &\text{if}\;\; \wh{\theta}\in p_h(D\setminus O_h); \\ P\left(A\right)&\text{otherwise}. \end{cases}
$$
Then $\{Q_{\wh{\theta}}\}_{\wh{\theta}\in\R}$ is a  disintegration of $P$ over $P_{\wh{\vT}}$ consistent with $\wh{\vT}$ and the following are all equivalent: 
\begin{enumerate}
\item
$N$ has the $P$-Markov property;

\item
$N$ has the $Q_{\wh{\theta}}$-Markov property for $P_{\wh{\vT}}\mbox{-a.a.}\;\;\wh{\theta}\in\R$;
\item
$N$ has the $Q_{\wh{\theta}}$-multinomial property for $P_{\wh{\vT}}\mbox{-a.a.}\;\;\wh{\theta}\in\R$.
\end{enumerate}
\end{cor}

{\bf{Proof.}} First note that by Lemma \ref{lem1} the family  $\{Q_{\wh{\theta}}\}_{\wh{\theta}\in\R}$ is a disintegration of $P$ over $P_{\wh{\vT}}$ consistent with $\wh{\vT}$.

Ad $(i)\Longrightarrow (ii)$: Assume that assertion $(i)$ holds true. Since $N$ is a eMRP$({\bf K}(h(\vT)))$, it follows by Theorem \ref{thm1} that $N$ is a MPP($\wh{\vT}$). Then according to \cite{lm1v}, Proposition 4.4, $N$ is a $Q_{\wh{\theta}}$-PP$(\wh{\theta})$ for $P_{\wh{\vT}}$-a.a. $\wh{\theta}\in\R$ and thus it has the $Q_{\wh{\theta}}$-Markov property.\smallskip

Ad $(ii)\Longrightarrow (i)$: Assume that assertion $(ii)$ holds true. Since $N$ is a  eMRP$({\bf K}(h(\vT)))$, it follows by Lemma \ref{lem3}  that $N$ is a MRP$({\bf{K}}(\wt{\vT}))$. Fix on an arbitrary $A\in\vS$ and put
$$
Q_{\wt{\theta}}\left(A\right):=\begin{cases}(P_{\bullet}(A)\circ h^{-1})(\wt{\theta}) &\text{if}\;\; \wt{\theta}\in V_h \\ P\left(A\right)& \text{otherwise,} \end{cases}
$$ 
where $\wt{\theta}:=h(\theta)$. By Lemma \ref{lem1} the family  $\{Q_{\wt{\theta}}\}_{\wt{\theta}\in\R}$ is a disintegration of $P$ over $P_{\wt{\vT}}$ consistent with $\wt{\vT}$. Applying now \cite{lm6z3}, Proposition 3.8, we obtain that there exists a $P_{\wt{\vT}}$-null set $U\in\B_k$ such that for any $\wt{\theta}\notin U$ the process $N$ is a $Q_{\wt{\theta}}$-RP$({\bf{K}}(\wt{\theta}))$. For any $\wh{\theta}\in\R$ put
$$
R_{\wh{\theta}}\left(A\right)=\begin{cases}Q_{\wt{\theta}}(A) &\text{if}\;\; \wh{\theta}:=p_k(\wt{\theta})\in U_{h,k}:=p_k(V_h\cup U^c) \\ P\left(A\right)& \text{otherwise.} \end{cases}
$$ 
Again by  Lemma \ref{lem1} the family  $\{R_{\wh{\theta}}\}_{\wh{\theta}\in\R}$ is a disintegration of $P$ over $P_{\wh{\vT}}$ consistent with $\wh{\vT}$, and for any $\wh{\theta}=p_k(\wt{\theta})\in U_{h,k}$   the process $W$ is $R_{\wh{\theta}}$-independent and ${\bf{K}}(\wh{\theta})=(R_{\wh{\theta}})_{W_n}
=(Q_{\wt{\theta}})_{W_n}={\bf{K}}(\wt{\theta})$. As a consequence we get that for any $\wh{\theta}\in U_{h,k}$ the process $N$ is a  $R_{\wh{\theta}}$-RP$({\bf{K}}(\wh{\theta}))$. But since $p_h(D\setminus O_h)$ is contained in $U_{h,k}$, we get $R_{\wh{\theta}}(A)=Q_{\wh{\theta}}(A)$ for any $\wh{\theta}\in p_h(D\setminus O_h)$. Thus, taking into account \cite{lm1v}, Proposition 4.4, we get that $N$ is a MPP$(\wh{\vT})$; hence $N$ has the $P$-Markov property (cf. e.g. \cite{Sc} Theorem 4.2.3).\smallskip

Ad $(ii)\Longrightarrow (iii)$: Assume that \textit{(ii)} holds true. Since \textit{(i)}$\Longleftrightarrow$\textit{(ii)}, it follows from Theorem \ref{thm1}, that \textit{(i)} is equivalent to the fact that $N$ is a MPP$(\wh{\vT})$. Applying now \cite{lm1v}, Proposition 4.4, we obtain that $N$ is a $Q_{\wh{\theta}}$-PP$(\wh{\theta})$  
for $P_{\wh{\vT}}$-a.a. $\wh{\theta}\in\R$. But the latter implies \textit{(iii)}.

 The implication $(iii)\Longrightarrow (ii)$ follows by an easy computation.
\hfill$\Box$ 

\begin{rems}\label{r2}
\normalfont
{\bf (a)} Assumption \ref{as}, more precisely its part concerning the differentiability of the distribution functions of $W_n$ with respect to $P_\theta$, is essential for the validity of Proposition \ref{prop1} and Theorem \ref{thm1}. \smallskip

In fact, consider  the trivial counting process $N$ defined by means of $N_{t}:= \left[t\right]$ for every $t\in\R_+$. It can be easily proven that for given $\theta_0>0$ the process $N$ is a Markov RP$(\mathbf{K}(\theta_0))$ with 
$$
{\bf{K}}(\theta_0)((-\infty,t]):=\begin{cases}0 &\text{if}\;\; t<1; \\ 1 & \text{if}\;\; t\geq 1, \end{cases}
$$
but not a Poisson process. 

For any $\theta\in\R$ define the set-function $P_\theta:\vS\longmapsto [0,1]$ by means of $P_{\theta}(A):=P(A)$ for any $A\in\vS$. It then can  easily be seen that the family $\{P_\theta\}_{\theta\in\R}$ is a disintegration of $P$ over $P_\vT$ consistent with any random variable $\vT$ such that $P_\vT(\{\theta_0\})=1$, and that $\{P_\theta\}_{\theta\in\R}$ and $N$ do not satisfy Assumption \ref{as}.\smallskip

{\bf (b)} A characterization of MPPs with mixing distribution $U$ (cf. e.g. \cite{SZ}, page 9 for the definition) in terms of the multinomial property has been obtained without any additional assumption by Schmidt and Zocher \cite{SZ}, Theorem 4.2. But it seems that such a characterization cannot be carried over to MPP$(\vT)$ without any additional assumption, since in general it is not possible given a distribution $U$ to find a real-valued random variable $\vT$ on $\vO$ with $P_\vT=U$, and a disintegration of $P$ over $U$ consistent with $\vT$ (compare Zocher \cite{Zo}, page 115). However, such an equivalence becomes possible under the essential assumptions of Proposition \ref{prop1} and Theorem \ref{thm1}. Thus, (b) together with (a) raises the following
\end{rems}

\begin{que}\label{que}
\normalfont
Is the assumption  of Proposition \ref{prop1} and Theorem \ref{thm1} concerning the existence of a disintegration of $P$ over $P_\vT$ consistent with $\vT$ necessary  for the validity of their conclusions ?
\end{que}

\section{Examples}\label{Examples} 
By $(\vO\times\vY,\vS\otimes{H},P\otimes{Q})$ we denote the product probability space of $(\vO,\vS,P)$ and $(\vY,H,Q)$, and by $\pi_{\vO}$ and $\pi_{\vY}$ the canonical projections from $\vO\times\vY$ onto $\vO$ and $\vY$, respectively.

{\em{Throughout what follows, we put $\varUpsilon:=(0,\infty)$, $H:=\mathfrak B(\varUpsilon)$, 
$\vO:=\varUpsilon^{\mathbb N}\times G$ for $G\in\B_d$, $\vS:=\mathfrak B(\vO)=\mathfrak B(\varUpsilon^{\mathbb N})\otimes \mathfrak B(G)$ for simplicity.}}\smallskip

First, we describe a method for the construction of non-trivial probability spaces admitting extended MRPs with mixing parameters $\vT$ and $h$, generalizing in this way Example 5.5 from \cite{lm6z3}.

\begin{ex}\label{ex2} 
\normalfont
Let  $\mu$ be an arbitrary probability measure on $\mathfrak B(G)$ and let $Q_{n}(\theta)$
be probability measures on $\mathfrak B(\varUpsilon)$ 
for all $n\in\mathbb N$ and for any fixed $\theta\in G$, which is absolutely continuous with respect to Lebesgue measure $\lambda$ on $\mathfrak B$. Suppose that there exists a measurable  map $h:G\longmapsto \R^k$ such that $Q_{n}(\theta)={\bf{K}}\left(h(\theta)\right)$ for any $n\in\N$, where for any $B\in \mathfrak B(\varUpsilon)$ the function ${\bf{K}}\left(h(\bullet)\right)(B):G\longmapsto\R$ is $\mathfrak B(G)$-measurable and ${\bf{K}}\left(h(\theta)\right)(\vY)=1$. It then follows that there exists a unique probability measure $\wt{P}_{\theta}:=\otimes_{n\in\mathbb N} Q_{n}(\theta)$ on $\mathfrak B(\varUpsilon^{\mathbb N})$. Put $P(E) :=\int  \wt{P}_{\theta}(E^{\theta})\mu(d\theta)$, for each $E\in\vS$, where $E^{\theta}$ is the $\theta$-section of $E$, and $P_{\theta}:=\wt{P}_{\theta}\otimes \delta_{\theta}$ for any $\theta\in G$, where $\delta_{\theta}$ is the Dirac measure at $\theta$. Then $P$ is a probability measure on $\vS$ and $\{\wt P_\theta\}_{\theta\in G}$ is a product regular conditional probability  on $\B(\vY^\N)$ (see \cite{smm}, Definition 1.1 for the definition and its properties); hence according to \cite{lm5} , Proposition 2.5, $\{P_{\theta}\}_{\theta\in G}$ is a disintegration of $P$ over $\mu$ consistent with $\pi_{G}$ (compare  \cite{lm6z3}, Example 5.5).

Clearly, putting $\vT:= \pi_{G}$ we get $P_{\vT}=\mu$. Set $W_{n}:= \pi_{n}$, where $\pi_n$ is the canonical projection from $\vO$ onto $\vY$, for any $n \in\mathbb N$ and $W:=\{W_{n}\}_{n\in\N}$.  It follows  
that $W$ is $P_\theta$-independent for any $\theta\in G$, implying together with \cite{lm1v}, Lemma 4.1, that $W$ is $P$-conditionally independent. Moreover, we have $(P_{\theta})_{W_{n}} =Q_{n}(\theta)={\bf K}(h(\theta))$ for all $n\in\mathbb N$ and $\theta\in G$, implying together with Lemma \ref{lem2} that for each $n\in\mathbb N$ the equality $P_{W_{n}\mid\vT}={\bf{K}}\left(h(\vT)\right)$ holds $P\upharpoonright\sigma(\vT)$-a.s. true. Put $T_n:=\sum_{k=1}^{n} W_k$ for any $n\in\N_0$ and $T:=\{T_n\}_{n\in\N_0}$, and let $N:=\{N_t\}_{t\in\R_+}$ be the counting process induced by $T$ by means of $N_t:=\sum_{n=1}^{\infty}\chi_{\{T_n\leq t\}}$ for all $t\in\R_+$ (cf. e.g \cite{Sc}, Theorem 2.1.1). Consequently, according to Definition \ref{rd} (a) the  counting process $N$  is a eMRP$({\bf K}(h(\vT)))$.
\end{ex}

In the next  examples it is shown that there exist non-trivial probability spaces satisfying all assumptions of  Theorem \ref{thm1}, which allow us to check whether  a  eMRP$({\bf K}(h(\vT)))$ is a Markov process or has the multinomial property.

\begin{exs}
\label{exs}
\normalfont
{\bf{(a)}} Let $G:=\vY\subseteq\R^d$ for d=1, let $\mu:={\bf{Ga}}(\alpha,\beta)$, with $\alpha,\beta>0$, be a probability measure on $\mathfrak B(\vY)$ and let $h:\vY\longmapsto\R$ be a function defined by means of $h(\theta):=a\theta+b$ for any $\theta>0$, where $a>0$ and $b\geq 0$ are constants. 
 Fix on arbitrary $\theta\in\vY$ and define the probability measures  $Q_n(\theta)$  on $\mathfrak B(\vY)$ by means of $Q_{n}(\theta):={\bf{Exp}}\left(h(\theta)\right)$ for all $n\in\N$. 
It then follows by Example \ref{ex2}, that there exist a map  $\vT:=\pi_{\vY}$, a probability measure $P$, a disintegration $\{P_\theta\}_{\theta\in\vY}$ of $P$ over $P_\vT=\mu$ consistent with $\vT$, and a counting process $N$ being a eMRP$({\bf K}(h(\vT)))$ such that its induced interarrival process $W$  satisfy condition $\left(P_\theta\right)_{W_n}=Q_{n}(\theta)$ for all $n\in\N$.

Define the map $C\in\mathcal L^1(P_{h(\vT)})$ by $C(h(\theta)):=h(\theta)$ for any $\theta\in\vY$, and for any fixed $\theta\in\vY$ define the density $f_{h(\theta)}:=F'_{h(\theta)}$  by $f_{h(\theta)}(t):=h(\theta)\cdot e^{-h(\theta)t}$ for any $t>0$. Clearly, for any fixed $\theta\in\vY$, the density $f_{h(\theta)}$ is dominated by $C(h(\theta))$, and the function $p_h:\vY\longmapsto \R$ defined by means of $p_h(\theta):=\lim_{t\rightarrow 0}f_{h(\theta)}(t)=h(\theta)$ for any $\theta\in\vY$, is positive and injective; hence $\{P_\theta\}_{\theta\in\vY}$, $N$ and $h$ satisfy all assumptions of Theorem \ref{thm1}.

Let $\wh{\vT}:=p_h\circ\vT$ and put $Q_{\wh{\theta}}(E):=\left(P_{\bullet}(E)\circ p^{-1}_h\right)(\wh{\theta})$ for any $\wh{\theta}>b$ and $E\in\vS$. Since $p_h=h$, we get $\wh\vT=\wt\vT$ and $Q_{\wh\theta}=Q_{\wt\theta}$ for any $\wh\theta=p_h(\theta)=h(\theta)=\wt\theta$, $\theta\in\vY$.
 
 By Lemma \ref{lem1} follows that  $\{Q_{\wh{\theta}}\}_{\wh{\theta}>b}$ is a disintegration of $P$ over $P_{\wh{\vT}}$ consistent with $\wh{\vT}$,  condition $\left(Q_{\wh{\theta}}\right)_{W_n}={\bf{Exp}}(\wh{\theta})$ holds true for any $n\in\N$ and $\wh{\theta}>b$, and the process $W$ is $Q_{\wh{\theta}}$-independent. But the latter implies that $N$ is $Q_{\wh{\theta}}$-PP$(\wh{\theta})$ for any $\wh{\theta}>b$ (cf. e.g. \cite{Sc}, Theorem 2.3.4). Thus, according to \cite{lm1v}, Proposition 4.4, we deduce that $N$ is a MPP$(\wh{\vT})$ implying together with Theorem \ref{thm1} that $N$ satisfies each of the equivalent conditions \textit{(i)}, \textit{(ii)} and \textit{(iii)} of Theorem \ref{thm1}.\smallskip

{\bf{(b)}}  Let $G$ and $\mu$ be as in (a) and let $h:\vY\longmapsto\R$ be a function defined by means of  $h(\theta):=\frac{1}{\theta}$ for any $\theta\in\vY$. Fix on an arbitrary $\theta\in\vY$ and take $Q_{n}(\theta):={\bf{Par}}\left(h(\theta),1\right)$ for all $n\in\N$, i.e.
$$
Q_{n}(\theta)(B):=\int_{B}\theta\cdot\left(\frac{1/\theta}{1/\theta+t}\right)^{2}\cdot\chi_{(0,\infty)}(t)\,\lambda(dt)\quad\text{for any}\; B\in\mathfrak B(\vY).
$$
It then follows by Example \ref{ex2}, that there exist a map  $\vT:=\pi_{\vY}$, a probability measure $P$, a disintegration $\{P_\theta\}_{\theta\in\vY}$ of $P$ over $P_\vT=\mu$ consistent with $\vT$, and a counting process $N$ being a eMRP$({\bf K}(h(\vT)))$ such that its induced interarrival process $W$  satisfy condition $\left(P_\theta\right)_{W_n}=Q_{n}(\theta)$ for all $n\in\N$.

Define the map $C\in\mathcal L^1(P_{h(\vT)})$ by $C(h(\theta)):=\theta$ for any $\theta\in\vY$, and for any fixed $\theta\in\vY$ define the density $f_{h(\theta)}:=F'_{h(\theta)}$  by $f_{h(\theta)}(t):=\theta\cdot\left(\frac{1/\theta}{1/\theta+t}\right)^{2}$ for any $t>0$. Clearly, for any fixed $\theta\in\vY$, the density $f_{h(\theta)}$ is dominated by $C(h(\theta))$, and the function $p_h:\vY\longmapsto \R$ defined by means of  $p_h(\theta):=\lim_{t\rightarrow 0}f_{h(\theta)}(t)=\theta$ for any $\theta\in\vY$, is positive and injective; hence $\{P_\theta\}_{\theta\in\vY}$, $N$ and $h$ satisfy all assumptions of Theorem \ref{thm1}.

Let $\wh{\vT}:=p_h\circ\vT$ and put $Q_{\wh{\theta}}(E):=\left(P_{\bullet}(E)\circ p^{-1}_h\right)(\wh{\theta})$ for any $\wh{\theta}>0$ and $E\in\vS$. Since $p_h= id_\vY$, we get $\wh\vT=\vT$ and $Q_{\wh\theta}=P_\theta$ for any $\wh\theta=\theta\in\vY$.

Assume, if possible, that $N$ is a Markov process. It then follows by Theorem \ref{thm1} that $N$ is a MPP$(\wh{\vT})$ or equivalently that $(Q_{\wh{\theta}})_{W_n}={\bf Exp}(\wh{\theta})$ for any $n\in\N$ and $W$ is $Q_{\wh{\theta}}$-independent for $P_{\wh{\vT}}$-a.a. $\wh{\theta}>0$ (see \cite{lm1v}, Proposition 4.5).
Moreover, since $P_\theta(E)=Q_{\wh{\theta}}(E)$ for all $E\in\vS$, we deduce that 
$$
{\bf Par}(h(\theta),1)=(P_\theta)_{W_n}=(Q_{\wh{\theta}})_{W_n}={\bf Exp}(\theta),
$$
a contradiction.
\end{exs}

It follows an example to show that the part of Huang's Assumption $(\ast)$ concerning the boundedness of $F'_{\wt y}$ by a constant $C>0$ is not necessary.

\begin{ex}
\label{ex4}
\normalfont
Let $G$, $\mu$ and $h$ be as in Example \ref{exs} (a). Since $N$ is an eMRP$({\bf Exp}(h(\vT)))$, it follows by  Lemma \ref{lem3} that it is a MRP$({\bf Exp}(\wt\vT))$. The latter together with \cite{lm6z3}, Theorem 4.9, yields that $N$ is a MRP$(\{Q_{\wt\theta}\},P_{\wt\vT})$.

According to Example \ref{exs} (a), $N$ and $\{P_\theta\}_{\theta\in\vY}$ satisfy all assumptions of Theorem \ref{thm1}; hence its conclusions. In particular $N$ is a   MPP$(\wh\vT)$; hence a MPP$(\wt\vT)$, because $p_h=h$. 

Assume that $N$ is a MPP$(\{Q_{\wt\theta}\},P_{\wt\vT})$. Then for $P_{\wt\vT}$-a.a. $\wt\theta>b$ we have that $(Q_{\wt\theta})_{W_n}={\bf Exp}(\wt\theta)$. But since Since $N$ is a MRP$({\bf Exp}(\wt\vT))$, it follows by \cite{lm1v}, Lemma 4.1, that $W$ is $Q_{\wt\theta}$-independent for $P_{\wt\vT}$-a.a. $\wt\theta>b$; hence$N$ is a $Q_{\wt\theta}$-PP$(\wt\theta)$ (cf. e.g. \cite{Sc}, Theorem 2.3.4), implying that $N$ has the $Q_{\wt\theta}$-Markov property (cf. e.g. \cite{Sc}, Corollary 3.1.2), equivalently $N$ has the $P$-Markov property, see Corollary \ref{cor1}. 

For the inverse implication of Theorem 3 from \cite{hu}, assume that $N$ has the $P$-Markov property. It the follows from Theorem \ref{thm1} that $N$ is a MPP$(\wt\vT)$; implying together with  Proposition \ref{prop22} that $N$ is a MPP$(\{Q_{\wt\theta}\},P_{\wt\vT})$. As a consequence, we get that the conclusions of \cite{hu}, Theorem 3 hold true.

But since for any $\wt\theta>b$ and $n\in\N$ we have $(Q_{\wt\theta})_{W_n}={\bf Exp}(\wt\theta)$, it follows that there does not exist any positive constant $C$ with $F'_{\wt\theta}(t)<C$ for all  $t>0$ and $\wt\theta>b$.

Thus the part of Assumption $(\ast)$ concerning the boundedness of $F'_{\wt y}$ by a constant $C>0$ is not necessary. In particular, in the case of the above example Huang's Theorem cannot be applied.
\end{ex}

\medskip

{\small
\noindent{\sc N.D. Macheras and S.M Tzaninis}\\
{\sc Department of Statistics and Insurance Science}\\
{\sc University of Piraeus, 80 Karaoli and Dimitriou street}\\
{\sc 185 34 Piraeus, Greece}\\
E-mail: {\tt macheras@unipi.gr}$\;$ {\sc and}$\;$ {\tt stzaninis@unipi.gr}
}

\end{document}